\documentclass{amsart}
\usepackage{amsmath}
\usepackage{amssymb}
\usepackage{amsthm}

\DeclareMathOperator{\Tr}{Tr}
\DeclareMathOperator{\tr}{tr}

\DeclareMathOperator{\Vol}{Vol}
\DeclareMathOperator{\dvol}{dvol}

\DeclareMathOperator{\Ric}{Ric}

\newcommand{\oJ}{\overline{J}}

\newcommand{\oP}{\overline{P}}

\newcommand{\of}{\overline{f}}
\newcommand{\og}{\overline{g}}

\newcommand{\odelta}{\overline{\delta}}

\newcommand{\oDelta}{\overline{\Delta}}

\newcommand{\oRic}{\overline{\Ric}}

\newcommand{\onabla}{\overline{\nabla}}

\newcommand{\cT}{\widetilde{T}}

\newcommand{\hf}{\widehat{f}}
\newcommand{\hg}{\widehat{g}}
\newcommand{\hh}{\widehat{h}}
\newcommand{\hu}{\widehat{u}}
\newcommand{\hA}{\widehat{A}}
\newcommand{\hB}{\widehat{B}}

\newcommand{\hH}{\widehat{H}}

\newcommand{\hL}{\widehat{L}}
\newcommand{\hP}{\widehat{P}}
\newcommand{\hQ}{\widehat{Q}}
\newcommand{\hT}{\widehat{T}}
\newcommand{\hmB}{\widehat{\mathcal{B}}}

\newcommand{\hmF}{\widehat{\mathcal{F}}}
\newcommand{\hmG}{\widehat{\mathcal{G}}}

\newcommand{\hmT}{\widehat{\mathcal{T}}}

\newcommand{\hnabla}{\widehat{\nabla}}

\newcommand{\heta}{\widehat{\eta}}

\newcommand{\hpsi}{\widehat{\psi}}

\newcommand{\lp}{\langle}
\newcommand{\rp}{\rangle}
\newcommand{\lv}{\lvert}
\newcommand{\rv}{\rvert}
\newcommand{\lV}{\lVert}
\newcommand{\rV}{\rVert}

\newcommand{\mB}{\mathcal{B}}
\newcommand{\mC}{\mathcal{C}}

\newcommand{\mE}{\mathcal{E}}
\newcommand{\mF}{\mathcal{F}}
\newcommand{\mG}{\mathcal{G}}

\newcommand{\mQ}{\mathcal{Q}}

\newcommand{\mT}{\mathcal{T}}

\newcommand{\bR}{\mathbb{R}}

\newcommand{\suchthat}{\mathrel{}\middle|\mathrel{}}

\def\sideremark#1{\ifvmode\leavevmode\fi\vadjust{\vbox to0pt{\vss
 \hbox to 0pt{\hskip\hsize\hskip1em
 \vbox{\hsize3cm\tiny\raggedright\pretolerance10000
 \noindent #1\hfill}\hss}\vbox to8pt{\vfil}\vss}}}

\newcommand{\comment}[1]{}

\newtheorem{thm}{Theorem}[section]
\newtheorem{prop}[thm]{Proposition}
\newtheorem{lem}[thm]{Lemma}
\newtheorem{cor}[thm]{Corollary}

\theoremstyle{definition}
\newtheorem{defn}[thm]{Definition}

\theoremstyle{remark}
\newtheorem{remark}[thm]{Remark}

\numberwithin{equation}{section}

\begin{document}

\title{Boundary operators associated to the Paneitz operator}
\author{Jeffrey S. Case}
\address{109 McAllister Building \\ Penn State University \\ University Park, PA 16801}
\email{jscase@psu.edu}
\begin{abstract}
We describe a set of conformally covariant boundary operators associated to the Paneitz operator, in the sense that they give rise to a conformally covariant energy functional for the Paneitz operator on a compact Riemannian manifold with boundary.  These operators naturally give rise to a first- and third-order conformally covariant pseudodifferential operator.  In the setting of Poincar\'e--Einstein manifolds, we show that these operators agree with the fractional GJMS operators of Graham and Zworski.  We also use our operators to establish some new sharp Sobolev trace inequalities.
\end{abstract}
\maketitle

\section{Introduction}
\label{sec:intro}

The Paneitz operator is a fourth-order conformally covariant differential operator which gives the simplest higher-order analogue of the conformal Laplacian.  For example, the conformal Laplacian and the Paneitz operator arise in the norm computations for the critical Sobolev embeddings $W^{1,2}(\bR^n)\hookrightarrow L^{\frac{2n}{n-2}}(\bR^n)$ and $W^{2,2}(\bR^n)\hookrightarrow L^{\frac{2n}{n-4}}(\bR^n)$, respectively, and in the study of the scalar curvature and the (fourth-order) $Q$-curvature, respectively, within a conformal class.

To study conformally covariant operators on a manifold with boundary requires specifying conformally covariant boundary operators.  Together, these operators are expected to control both critical Sobolev trace embeddings and certain scalar invariants.  This is completely understood for the conformal Laplacian, and it is the purpose of this article to completely describe the analogous properties for the Paneitz operator.

Given a Riemannian manifold $(X^{n+1},g)$ with boundary $(M^n,h)=(\partial X,g\rv_{TM})$, the conformal Laplacian $L_2$, the trace $B_0^1$, and the conformal Robin operator $B_1^1$ are defined by
\[ L_2u = -\Delta u + \frac{n-1}{2}Ju, \qquad B_0^1u = u\rv_M, \qquad B_1^1u = \eta u + \frac{n-1}{2n}Hu\rv_M \]
for all $u\in C^\infty(X)$, where $J$ is the trace of the Schouten tensor $P=\frac{1}{n-1}(\Ric-Jg)$, $\eta$ is the outward-pointing normal, and $H=\tr_h\nabla\eta$ is the mean curvature of $M$.  It follows that
\[ \int_X u_1L_2(u_2) + \oint_M f_1\,B_1^1(u_2) = \int_X \left(\lp\nabla u_1,\nabla u_2\rp + \frac{n-1}{2}Ju_1u_2\right) + \frac{n-1}{2n}\oint_M Hf_1f_2 \]
for all $u_1,u_2\in C^\infty(X)$, where $f_j=B_0^1(u_j)$ for $j\in\{1,2\}$.  Denote this pairing by $\mQ_2$.  If $\hg=e^{2\sigma}g$ is a conformal rescaling of $g$, then
\begin{equation}
 \label{eqn:l2_conf}
 \begin{split}
  \hL_2u & = e^{-\frac{n+3}{2}\sigma}L_2\left(e^{\frac{n-1}{2}\sigma}u\right), \\
  B_k^1u & = e^{-\frac{n+2k-1}{2}\sigma}B_k^1\left(e^{\frac{n-1}{2}\sigma}u\right)  
 \end{split}
\end{equation}
for all $u\in C^\infty(X)$ and all $k\in\{0,1\}$.  It follows that $\mQ_2$ is a conformally covariant symmetric bilinear form.  Let $\mE_2(u)=\mQ_2(u,u)$ be the associated energy.  Then the boundary Yamabe invariant
\[ Y_2(X,M) = \inf \left\{ \mE_2(u) \suchthat u\in C^\infty(X), \oint_M \left\lv B_0^1(u)\right\rv^{\frac{2n}{n-1}}=1 \right\} \]
is a conformal invariant.

A number of analytic and geometric consequences can be derived from $Y_2(X,M)$.  Using conformal covariance and an Obata-type theorem, Escobar~\cite{Escobar1988} computed $Y_2(\bR_+^{n+1},\bR^n)$, thus giving a norm computation for the Sobolev trace embedding $W^{1,2}(\bR_+^{n+1})\hookrightarrow L^{\frac{2n}{n-1}}(\bR^n)$.  Positive smooth critical points of $\mE_2$ subject to the constraint defining $Y_2(X,M)$ give rise to metrics which are scalar flat in $X$ and make $M$ have constant mean curvature~\cite{Escobar1988}.  In fact, when $Y_2(X,M)$ is finite, it can be realized by a positive smooth function~\cite{Almaraz2010,Escobar1992a,Marques2005,Marques2007}.

Suppose that $X$ is compact and that the bottom of the Dirichlet spectrum of the conformal Laplacian,
\[ \lambda_1(L_2) = \inf\left\{ \mE_2(u) \suchthat u\in C^\infty(X), B_0^1u=0, \int_X u^2=1 \right\} , \]
is positive; equivalently, suppose that $Y_2(X,M)$ is finite~\cite{Escobar1992ac}.  Then the energy $\mE_2$ is uniquely minimized in the class of functions $u\in C^\infty(X)$ with fixed trace $f=B_0^1u$ by a function $u\in\ker L_2$.  This gives a norm computation for the Sobolev trace embedding theorem: The map $C^\infty(X)\ni u\mapsto u\rv_M$ extends to a surjective bounded linear mapping $W^{1,2}(X)\to H^{1/2}(M)$ (cf.\ \cite{Treves1980}).

One can also use the pair $(L_2,B_1^1)$ to define a first order conformally covariant pseudodifferential operator as follows (cf.\ \cite{Branson1997}): Suppose $\ker L_2\cap\ker B_0^1=\{0\}$.  Then to each $f\in C^\infty(X)$, we define $\mB_1^1f:=B_1^1u_f$ for $u_f\in C^\infty(X)$ the unique function such that $L_2u_f=0$ and $B_0^1u=f$.  It follows from~\eqref{eqn:l2_conf} that $\mB_1^1$ is conformally covariant; indeed,
\[ \hmB_1^1f = e^{-\frac{n+1}{2}\sigma}\mB_1^1\left(e^{\frac{n-1}{2}\sigma}f\right) . \]
One can also define a first order conformally covariant pseudodifferential operator on $M$, namely the fractional GJMS operator $P_1$, via scattering theory~\cite{GrahamZworski2003}.  These operators are the same~\cite{GuillarmouGuillope2007}.  In particular, $\mB_1^1=(-\Delta)^{1/2}$ on Euclidean space~\cite{Branson1997,CaffarelliSilvestre2007,ChangGonzalez2011}.

Consider now the Paneitz operator $L_4\colon C^\infty(X)\to C^\infty(X)$ defined by
\[ L_4u = \Delta^2u + \delta\left((4P-(n-1)Jg)(\nabla u)\right) + \frac{n-3}{2}Q_4u \]
for all $u\in C^\infty(X)$, where $\delta=\tr\nabla$ is the negative of the formal adjoint of the Levi-Civita connection and
\[ Q_4 = -\Delta J - 2\lv P\rv^2 + \frac{n+1}{2}J^2 \]
is the (fourth-order) $Q$-curvature.  The Paneitz operator is conformally covariant~\cite{Paneitz1983}; indeed, if $\hg=e^{2\sigma}g$, then
\begin{equation}
 \label{eqn:l4_conf}
 \hL_4(u) = e^{-\frac{n+5}{2}\sigma}L_4\left(e^{\frac{n-3}{2}\sigma}u\right) .
\end{equation}
As a fourth-order operator, the Paneitz operator should have four associated boundary operators.  Define operators $B_k^3\colon C^\infty(X)\to C^\infty(M)$ for $k\in\{0,1,2,3\}$ by
\begin{align*}
 B_0^3u & = u , \\
 B_1^3u & = \eta u + \frac{n-3}{2n}Hu, \\
 B_2^3u & = -\oDelta u + \nabla^2u(\eta,\eta) + \frac{n-2}{n}H\eta u + \frac{n-3}{2}T_2^3u, \\
 B_3^3u & = -\eta\Delta u - 2\oDelta\eta u - \frac{n-3}{2n}H\nabla^2u(\eta,\eta) + \frac{4}{n-1}\lp A_0,\onabla^2u\rp - \frac{3n-5}{2n}H\oDelta u \\
  & \quad - \frac{2(n-4)}{n}\lp\onabla H,\onabla u\rp + S_2^3\eta u + \frac{n-3}{2}T_3^3u
\end{align*}
for all $u\in C^\infty(X)$, where $A=\nabla\eta\rv_{TM}$ is the second fundamental form, $A_0=A-\frac{H}{n}h$ is its tracefree part,
\begin{align*}
 S_2^3 & = -\frac{3n^2-7n+6}{4n^2}H^2 + \frac{n-7}{2}P(\eta,\eta) + \frac{3n-5}{2}\oJ + \frac{1}{2}\lv A_0\rv^2 , \\
 T_2^3 & = \oJ - P(\eta,\eta) + \frac{n-2}{2n^2}H^2, \\
 T_3^3 & = \eta J - \frac{2}{n}\oDelta H - \frac{4}{n-1}\lp A_0,\oP\rp + \frac{n-3}{2n}HP(\eta,\eta) \\
  & \quad + \frac{3n-1}{2n}H\oJ + \frac{n+1}{2n(n-1)}H\lv A_0\rv^2 - \frac{n^2-n+2}{4n^3}H^3 ,
\end{align*}
and barred operators are defined with respect to $(M^n,h)$.  These are boundary operators associated to the Paneitz operator in the following sense:

\begin{thm}
 \label{thm:boundary}
 Let $(X^{n+1},g)$ be a Riemannian manifold with boundary $M^n=\partial X$ and $n\geq2$.  Given $\sigma\in C^\infty(X)$, set $\hg=e^{2\sigma}g$.
 \begin{enumerate}
  \item Given any $u\in C^\infty(X)$ and any $k\in\{0,1,2,3\}$, it holds that
  \[ \hB_k^3(u) = e^{-\frac{n+2k-3}{2}\sigma\rv_M}B_k^3\left(e^{\frac{n-3}{2}\sigma}u\right) . \]
  \item Define $\mQ_4\colon C^\infty(X)\times C^\infty(X)\to\bR$ by
  \[ \mQ_4(u,v) = \int_X u\,L_4v + \oint_M \left( B_0^3(u)B_3^3(v) + B_1^3(u)B_2^3(v)\right) \]
  for all $u,v\in C^\infty(X)$.  Then $\mQ_4$ is a symmetric bilinear form.
 \end{enumerate}
\end{thm}

That $\mQ_4$ is a symmetric bilinear form follows from the identity
\begin{equation}
 \label{eqn:l4_symmetric}
 \begin{split}
  \mQ_4(u,v) & = \int_X \left((\Delta u)(\Delta v) - (4P-(n-1)Jg)(\nabla u,\nabla v) + \frac{n-3}{2}Q_4uv\right) \\
  & \quad + \oint_M\biggl( 2\lp\onabla f_1,\onabla\psi_2\rp + 2\lp\onabla f_2,\onabla\psi_1\rp - \frac{2}{n}H\psi_1\psi_2 \\
   & \qquad + \frac{n-3}{2}\cT_2^3\left(f_1\psi_2+f_2\psi_1\right) - \left(\frac{4}{n-1}A_0-\frac{2}{n}Hh\right)(\onabla f_1,\onabla f_2) \\
   & \qquad + \frac{n-3}{2}\left(\frac{2}{n}\oDelta H + T_3^3 - \frac{n-3}{2n}H\cT_2^3\right)f_1f_2\biggr)
 \end{split}
\end{equation}
for all $u,v\in C^\infty(X)$, where $f_1=B_0^3(u)$, $f_2=B_0^3(v)$, $\psi_1=B_1^3(u)$, $\psi_2=B_1^3(v)$, and
\[ \cT_2^3=T_2^3+\frac{2}{n^2}H^2 . \]
Theorem~\ref{thm:boundary} and the identity~\eqref{eqn:l4_symmetric} together justify our description of the operators $B_k^3$ as the boundary operators associated to the Paneitz operator.  However, this does not uniquely determine the boundary operators; see Remark~\ref{rk:uniqueness}.

Our main contribution in Theorem~\ref{thm:boundary} is the formula~\eqref{eqn:l4_symmetric}.  Motivated by questions involving the functional determinant, Chang and Qing~\cite{ChangQing1997a} studied boundary operators for the Paneitz operator on four-manifolds with boundary when acting on functions whose normal derivative vanishes along the boundary, and identified the operator $B_3^3$ and the identity~\eqref{eqn:l4_symmetric} in this setting (cf.\ Section~\ref{sec:3d}).  Branson and Gover~\cite{BransonGover2001} gave formulas for the boundary operators in the noncritical dimensions $n\geq4$ (cf.\ \cite{Gover2007}).  Grant~\cite{Grant2003} extended the work of Branson and Gover to the critical dimension four for the Paneitz operator, and Juhl~\cite{Juhl2009} verified the top order boundary operator in this case is $B_3^3$.  Grant also identified a local formula for $B_2^3$.  Stafford~\cite{Stafford2006} (see also~\cite{Juhl2009}) used the tractor expression of Branson and Gover to give the local formula for $B_3^3$ in higher dimensions.  A refinement of the tractor approach has been developed by Gover and Peterson~\cite{GoverPeterson2015}.

Theorem~\ref{thm:boundary} and~\eqref{eqn:l4_conf} implies that the energy $\mE_4(u):=\mQ_4(u,u)$ is a conformally covariant quadratic form.  In particular,
\begin{align*}
 Y_{4,1}(X,M) & = \inf \left\{ \mE_4(u) \suchthat u\in C^\infty(X), B_1^3u=0, \oint_M \left\lv B_0^3u\right\rv^{\frac{2n}{n-3}}=1 \right\}, \\
 Y_{4,2}(X,M) & = \inf \left\{ \mE_4(u) \suchthat u\in C^\infty(X), \oint_M \left\lv B_0^3u\right\rv^{\frac{2n}{n-3}}=1 \right\}
\end{align*}
are conformal invariants when $n\geq 4$.  Extremizing $Y_{4,1}(X,M)$ or $Y_{4,2}(X,M)$ provides one approach to constructing conformal metrics while prescribing the interior $Q$-curvature and certain combinations of the mean curvature and the curvatures $T_2^3$ and $T_3^3$ along the boundary.

\begin{prop}
 \label{prop:l4_geometric_euler}
 Let $(X^{n+1},g)$ be a Riemannian manifold with boundary $M^n=\partial X$ and $n\geq4$.  Suppose that $u\in C^\infty(X)$ is a positive admissible function such that $\mE_4(u)=Y_{4,k}(X,M)$ for $k\in\{1,2\}$.  Then the metric $\hg=u^{\frac{4}{n-3}}g$ is such that $\hQ_4=0$, $\hT_k^3=0$, and $\hT_3^3=\frac{2}{n-3}Y_{4,k}(X,M)$.
\end{prop}

Here $T_1^3=\frac{1}{n}H$ and $u$ is admissible if it is an admissible test function for the definition of $Y_{4,k}(X,M)$.  Proposition~\ref{prop:l4_geometric_euler} follows easily from Theorem~\ref{thm:boundary}, and in a similar way one can obtain a variational characterization for some other conformal curvature prescription problems.

While we do not consider here the general problems of extremizing $Y_{4,1}(X,M)$ or $Y_{4,2}(X,M)$, we do compute its value on the round upper hemisphere, or equivalently, Euclidean upper half space.  The first step in this direction is to relate the energy $\mE_4(u)$ of a Paneitz-harmonic function to certain boundary energies.  This can be done whenever the bottom of the Dirichlet spectrum of the Paneitz operator,
\[ \lambda_1(L_4) = \inf\left\{ \mE_4(u) \suchthat u\in C^\infty(X), B_0^3(u)=0=B_1^3(u), \int_X u^2=1 \right\} , \]
is positive.  Note that this is a conformally invariant condition.  The desired relation comes from the following extension result.

\begin{prop}
 \label{prop:extension}
 Let $(X^{n+1},g)$ be a compact Riemannian manifold with boundary $M^n=\partial X$ and $n\geq4$, and suppose that $\lambda_1(L_4)>0$.  Then for each $f,\psi\in C^\infty(M)$, there exists a unique function $u_{f,\psi}\in C^\infty(X)$ such that
 \begin{equation}
  \label{eqn:extension}
  \begin{cases}
   L_4u_{f,\psi} = 0, \\
   B_0^3u_{f,\psi} = f, \\
   B_1^3u_{f,\psi} = \psi .
  \end{cases}
 \end{equation}
 Indeed, $u_{f,\psi}$ is characterized by the property
 \begin{equation}
  \label{eqn:extension_energy}
  \mE_4(u_{f,\psi}) = \inf\left\{ \mE_4(u) \suchthat u\in C^\infty(X), B_0^3u=f, B_1^3u=\psi \right\} .
 \end{equation}
\end{prop}

A key point is that, by Theorem~\ref{thm:boundary}, the boundary value problem $(L_4;(B_0^3,B_1^3))$ defined by~\eqref{eqn:extension} is formally self-adjoint and satisfies the Lopatinskii--Shapiro conditions (cf.\ \cite{AgmonDouglisNirenberg1964,BransonGover2001}).  Hence weak solutions of~\eqref{eqn:extension} are always smooth, while the spectral assumption guarantees that the infimum~\eqref{eqn:extension_energy} is finite.

Suppose more generally that the $L^2$-kernel of the Paneitz operator on $(X^{n+1},g)$ is trivial; i.e.\ $\ker L_4\cap\ker B_0^3\cap\ker B_1^3=\{0\}$.  Given $f,\psi\in C^\infty(M)$, there are unique functions $u_{f,0},u_{0,\psi}\in C^\infty(X)$ such that~\eqref{eqn:extension} holds.  Therefore the operators $\mB_1^3,\mB_3^3\colon C^\infty(M)\to C^\infty(M)$ defined by $\mB_1^3\psi=\frac{1}{2}B_2^3u_{0,\psi}$ and $\mB_3f=\frac{1}{2}B_3^3u_{f,0}$ are well-defined.  Indeed, $\mB_1^3$ and $\mB_3^3$ are formally self-adjoint conformally covariant first- and third-order pseudodifferential operators with principal symbol that of $(-\Delta)^{1/2}$ and $(-\Delta)^{3/2}$, respectively; see Section~\ref{sec:asymptotics}.  On the other hand, if $g$ is conformal to a complete Einstein metric $g_+$ with $\Ric(g_+)=-ng_+$ in the interior of $X$ --- briefly, if $(X^{n+1},M^n,g_+)$ is a Poincar\'e--Einstein manifold --- then scattering theory defines first- and third-order pseudodifferential operators $P_1$ and $P_3$ with the same properties as $\mB_1^3$ and $\mB_3^3$, respectively~\cite{GrahamZworski2003}.  These operators are the same:

\begin{thm}
 \label{thm:asymptotics}
 Let $(X^{n+1},M^n,g_+)$ be a Poincar\'e--Einstein manifold with $n\geq2$ and $\frac{n^2-1}{4},\frac{n^2-9}{4}\not\in\sigma_{pp}(-\Delta_{g_+})$.  Let $\rho$ be a defining function for $M$ and set $g=\rho^2g_+$.  Let $u\in C^\infty(X)$ be such that $L_4u=0$.  Set $f=B_0^3u$ and $\psi=B_1^3u$.  Then
 \begin{align*}
  B_3^3u & = 2P_3f, \\
  B_2^3u & = 2P_1\psi .
 \end{align*}
 In particular, $\mB_3^3=P_3$ and $\mB_1^3=P_1$.
\end{thm}

A surprising consequence of Theorem~\ref{thm:asymptotics} is that, at least for Poincar\'e--Einstein manifolds, one can compute $P_3f$ without first finding the unique extension $u_{f,0}$; one need only find some Paneitz-harmonic function $u$ with $B_0^3u=f$ (cf.\ \cite{CaseChang2013,ChangGonzalez2011}).  This generalizes the observation that this property holds for hyperbolic space~\cite{AcheChang2015,BransonGover2001}.

Though both the pairs $\mB_1^3,\mB_3^3$ and $P_1,P_3$ are defined for all manifolds with boundary, it is only in the Poincar\'e--Einstein case that we can verify that they coincide.  The proof of Theorem~\ref{thm:asymptotics} relies on the observation that, when computed with respect to $g_+$, the Paneitz operator factors through the operators used in the Poisson equations defining $P_1$ and $P_3$ via scattering,
\begin{equation}
 \label{eqn:paneitz_einstein}
 \left(L_4\right)_{g_+} = \left(-\Delta_{g_+} - \frac{n^2-1}{4}\right)\circ\left(-\Delta_{g_+} - \frac{n^2-9}{4}\right) .
\end{equation}
If the Paneitz operator factors as in~\eqref{eqn:paneitz_einstein} for some metric $g$, then $g$ is Einstein~\cite{Robert2009}.

Combining Proposition~\ref{prop:extension} and Theorem~\ref{thm:asymptotics} yields the following relationship between the Paneitz energy of a function in $X$ and the energy of its trace:

\begin{cor}
 \label{cor:l4_sobolev_trace}
 Let $(X^{n+1},M^n,g_+)$ be a Poincar\'e--Einstein manifold such that $\lambda_1(L_4)>0$ and $n\geq2$.  Let $\rho$ be a defining function for $M$, set $g=\rho^2g_+$, and fix $f,\psi\in C^\infty(M)$.  Then
 \[  \frac{1}{2}\mE_4(u) \geq \oint_M f\,P_3f + \oint_M \psi\,P_1\psi \]
 for all $u\in C^\infty(X)$ with $B_0^3u=f$ and $B_1^3u=\psi$.  Moreover, equality holds if and only if $u=u_{f,\psi}$ for $u_{f,\psi}$ as in~\eqref{eqn:extension}.
\end{cor}

Recall the Sobolev trace theorem: In Euclidean upper half space, the map $u\mapsto (B_0^3u,B_1^3u)$ can be extended to a bounded linear mapping
\[ \Tr \colon W^{2,2}(\bR_+^{n+1}) \to H^{3/2}(\bR^n) \oplus H^{1/2}(\bR^n) \]
which admits a bounded right inverse (cf.\ \cite{Treves1980}).  Corollary~\ref{cor:l4_sobolev_trace} and Reilly's formula~\cite{Reilly1977} together give a norm computation for this result (cf.\ Corollary~\ref{cor:trace}).

Using the sharp fractional Sobolev inequalities~\cite{Beckner1993,Lieb1983}, we obtain from Corollary~\ref{cor:l4_sobolev_trace} a sharp Sobolev trace inequality which gives a norm computation for the embedding $W^{2,2}(\bR_+^{n+1})\hookrightarrow L^{\frac{2n}{n-3}}(\bR^n)\oplus L^{\frac{2n}{n-1}}(\bR^n)$.  In particular, this result computes the constants $Y_{4,1}(\bR_+^{n+1},\bR^n)$ and $Y_{4,2}(\bR_+^{n+1},\bR^n)$ and also characterizes their extremal functions:

\begin{thm}
 \label{thm:sobolev}
 Let $n\geq4$ and let $(S_+^{n+1},d\theta^2)$ be the upper hemisphere
 \[ S_+^{n+1} = \left\{ x=(x_0,\dotsc,x_{n+1})\in\bR^{n+2} \suchthat \lv x\rv^2=1, x_0\geq0 \right\} \]
 with the round metric $d\theta^2=dx^2\rv_{TS_+^{n+1}}$ of constant sectional curvature one.  Given any $u\in C^\infty(S_+^{n+1})$, it holds that
 \begin{multline*}
  \int_{S_+^{n+1}} \left((\Delta u)^2 + \frac{n^2-5}{2}\lv\nabla u\rv^2 + \frac{(n^2-1)(n^2-9)}{16}u^2\right) \\ \geq 4\oint_{S^n} \left(\psi\oDelta f - \frac{(n-1)(n-3)}{8}f\psi\right) + C_3\lV f\rV_{\frac{2n}{n-3}}^2 + C_1\lV\psi\rV_{\frac{2n}{n-1}}^2 ,
 \end{multline*}
 where $f=B_0^3u$, $\psi=B_1^3u$, and
 \[ C_k = 2^{k+1}\pi^{\frac{k}{2}}\frac{\Gamma\left(\frac{n+k}{2}\right)}{\Gamma\left(\frac{n-k}{2}\right)}\left(\frac{\Gamma(n/2)}{\Gamma(n)}\right)^{\frac{k}{n}} \]
 for $k\in\{1,3\}$.  Equality holds if and only if $L_4u=0$, $f(x)=a_1\left(1+\xi_1\cdot x\right)^{-\frac{n-3}{2}}$ and $\psi(x)=a_2\left(1+\xi_2\cdot x\right)^{-\frac{n-1}{2}}$ for constants $a_1,a_2\in\bR$ and points $\xi_1,\xi_2\in B_1(0)\subset\bR^{n+1}$.  In particular,
 \[ Y_{4,1}\left(\bR_+^{n+1},\bR^n\right) = Y_{4,2}\left(\bR_+^{n+1},\bR^n\right) = C_3 . \]
\end{thm}

More generally, if $g$ is any conformally flat metric on $S_+^{n+1}$, it holds that
\begin{equation}
 \label{eqn:l4_sobolev_general}
 \mE_4(u) \geq C_3\left(\oint_{S^n} \lv f\rv^{\frac{2n}{n-3}}\right)^{\frac{n-3}{n}} + C_1\left(\oint_{S^n} \lv\psi\rv^{\frac{2n}{n-1}}\right)^{\frac{n-1}{n}}
\end{equation}
for all $u\in C^\infty(S_+^{n+1})$, where $f=B_0^3u$ and $\psi=B_1^3u$.  Restricting~\eqref{eqn:l4_sobolev_general} to functions $u\in C^\infty(S_+^{n+1})$ with $B_1^3u=0$ recovers known sharp Sobolev trace inequalities: Combining results of Lieb~\cite{Lieb1983} and R.\ Yang~\cite{Yang2013} gives a sharp Sobolev trace inequality on Euclidean upper half space; Ache and Chang~\cite{AcheChang2015} established analogous inequalities in the flat ball $B_1(0)\subset\bR^{n+1}$ and in the upper hemisphere with both the round metric and the adapted metric~\cite{CaseChang2013}.

The restriction $n\geq4$ in Proposition~\ref{prop:l4_geometric_euler} and Theorem~\ref{thm:sobolev} is somewhat artificial.  There is a natural modification of the energy $\mE_4$ defined on four-manifolds with boundary which incorporates the boundary operators $B_k^3$ for $k\in\{0,1,2,3\}$.  This functional is related to log determinant formulae~\cite{BransonGilkey1994,ChangQing1997a} and gives rise to analogues of $Y_{4,1}$ and $Y_{4,2}$.  In particular, it extends Proposition~\ref{prop:l4_geometric_euler} and Theorem~\ref{thm:sobolev} to this setting; for details, see Section~\ref{sec:3d}.

This article is organized as follows:

In Section~\ref{sec:commute} we recall some useful facts from conformal geometry and the study of manifolds with boundary.

In Section~\ref{sec:operators} we prove Theorem~\ref{thm:boundary}, derive the Euler equation for $\mE_4$, and thereby establish Proposition~\ref{prop:l4_geometric_euler}.

In Section~\ref{sec:asymptotics} we prove Proposition~\ref{prop:extension} and discuss the properties of the operators $\mB_1^3$ and $\mB_3^3$.  In particular, we prove Theorem~\ref{thm:asymptotics}.

In Section~\ref{sec:sobolev} we describe the relationship to sharp Sobolev trace inequalities, and in particular prove Theorem~\ref{thm:sobolev}.

In Section~\ref{sec:3d} we describe the particulars of the case of compact four-manifolds with boundary.

\subsection*{Acknowledgments}

I would like to thank Antonio Ache and Alice Chang for many discussions related to their article~\cite{AcheChang2015}, and also Rod Gover for bringing my attention to the theses~\cite{Grant2003,Stafford2006}.

\section{Manifolds with boundary}
\label{sec:commute}

Throughout this article, we exclusively study manifolds $(X^{n+1},g)$ with boundary $(M^n,h)=(\partial X,g\rv_{TM})$.  We always denote intrinsic geometric quantities on $(M^n,h)$ with bars; e.g.\ $\oP=\frac{1}{n-2}(\oRic - \oJ h)$ is the Schouten tensor of $(M^n,h)$ and $\oJ=\tr_h\oP$ is its trace.  Integrals over $X$ with respect to the volume element induced by $g$ are denoted by $\int$; integrals over $M$ with respect to the volume element induced by $h$ are denoted by $\oint$.  When the context is clear, we omit restriction symbols (e.g.\ Proposition~\ref{prop:oi01}).

When studying boundary operators for the Paneitz operator, we often need to relate the intrinsic and extrinsic geometry of the boundary $M$ of $X$.  The first such relationship comes from the Gauss--Codazzi equations, stated here in terms of the Schouten tensor.  Since these equations are well-known, we omit their proof.

\begin{lem}
 \label{lem:gauss_codazzi}
 Let $(X^{n+1},g)$ be a Riemannian manifold with boundary $M^n=\partial X$ and let $\eta$ denote the outward-pointing normal along $M$.  Then
 \begin{align*}
  J & = \oJ + P(\eta,\eta) - \frac{1}{2n}H^2 + \frac{1}{2(n-1)}\lv A_0\rv^2 , \\
  P & = \oP - \frac{1}{n}HA_0 - \frac{1}{2n^2}H^2h + \frac{1}{n-2}\left(W(\eta,\cdot,\eta,\cdot) + A_0^2 - \frac{1}{2(n-1)}\lv A_0\rv^2 h\right) .
 \end{align*}
 Moreover, as sections of $T^\ast M$,
 \[ \odelta A = (n-1)P(\eta,\cdot) + dH . \]
\end{lem}

The other relationships we need involve second- and third-order derivatives of functions $u\in C^\infty(X)$ when evaluated on $M$.  The relationship between the Laplacians of $(X^{n+1},g)$ and $(M^n,h)$ is well-known.  Again, we omit the proof.

\begin{lem}
 \label{lem:2derivatives}
 Let $(X^{n+1},g)$ be a Riemannian manifold with boundary $M^n=\partial X$ and let $\eta$ denote the outward-pointing normal along $M$.  Given any $u\in C^\infty(X)$ it holds that
 \[ \Delta u = \oDelta u + \nabla^2u(\eta,\eta) + H\eta u . \]
\end{lem}

Given $u,v\in C^\infty(X)$, the derivative $\eta\lp\nabla u,\nabla v\rp$ is decomposed into horizontal and normal components via the following lemma.

\begin{lem}
 \label{lem:derivatives_eta}
 Let $(X^{n+1},g)$ be a Riemannian manifold with boundary $M^n=\partial X$.  Given any $u,v\in C^\infty(X)$ it holds that
 \[ \nabla^2u(\eta,\nabla v) = \lp\onabla\eta u,\onabla v\rp + (\eta v)\nabla^2u(\eta,\eta) - A(\onabla u,\onabla v) . \]
\end{lem}

\begin{proof}
 We compute that
 \begin{align*}
  \nabla^2u(\eta,\nabla v) & = \lp\nabla\eta u,\nabla v\rp - \lp\nabla_{\nabla v}\eta,\nabla u\rp \\
  & = \lp\onabla\eta u,\onabla v\rp + (\eta v)\left(\eta\eta u - (\nabla_\eta\eta)u\right) - \lp\nabla_{\onabla v}\eta,\nabla u\rp . \qedhere
 \end{align*}
\end{proof}

Finally, we require the relationship between the normal derivative of the interior Laplacian $\eta\Delta u$ and the third-order normal derivative $\nabla^3u(\eta,\eta,\eta)$ for $u\in C^\infty(X)$.

\begin{lem}
 \label{lem:3derivatives}
 Let $(X^{n+1},g)$ be a Riemannian manifold with boundary $M^n=\partial X$.  Given any $f\in C^\infty(X)$, it holds that
 \begin{multline*}
  \eta\Delta f = \nabla^3f(\eta,\eta,\eta) + H\nabla^2f(\eta,\eta) + \oDelta\eta f - 2\odelta\left(A(\onabla f)\right) + \lp\onabla H,\onabla f\rp \\ - \left(\oJ + nP(\eta,\eta) + \frac{1}{2n}H^2 + \frac{2n-1}{2(n-1)}\lv A_0\rv^2\right)\eta f .
 \end{multline*}
\end{lem}

\begin{proof}
 Let $p\in M$ and let $\{e_i\}_{i=1}^n$ be an orthonormal basis for $T_pM$, extended to a neighborhood of $p$ in $X$ by parallel translation.  It follows that, at $p$,
 \begin{align*}
  \eta\Delta f & = \nabla^3f(\eta,\eta,\eta) - \Ric(\eta,\nabla f) + \sum_{i=1}^n\nabla^3f(e_i,e_i,\eta) , \\
  \odelta A & = -H\nabla_\eta\eta + \sum_{i=1}^n\onabla_{e_i}\nabla_{e_i}\eta , \\
  \lp A, \onabla^2f\rp & = \sum_{i=1}^n \left[ e_i\left(\nabla_{e_i}\eta\right)f - \left(\onabla_{e_i}\nabla_{e_i}\eta\right)f \right] .
 \end{align*}
 We then compute that
 \begin{align*}
  \sum_{i=1}^n\nabla^3f(e_i,e_i,\eta) & = \sum_{i=1}^n \left[ e_ie_i\eta f - 2e_i\left(\nabla_{e_i}\eta\right) + \left(\nabla_{e_i}\nabla_{e_i}\eta\right)f \right] \\
  & = \oDelta\eta f + H\nabla^2f(\eta,\eta) - 2\lp A,\onabla^2f\rp - \lp\odelta A,\onabla f\rp - \lv A\rv^2\eta f .
 \end{align*}
 The result is now an immediate consequence of Lemma~\ref{lem:derivatives_eta}.
\end{proof}
\section{Boundary operators for the Paneitz operator}
\label{sec:operators}

The claim that the operators $B_k^3$ for $k\in\{0,1,2,3\}$ are boundary operators associated to the weighted Paneitz operator includes two statements.  First, the boundary operators $B_k^3$ are conformally covariant and act on densities of the same weight as the Paneitz operator; i.e.\ the precomposition factor is the same as in~\eqref{eqn:l4_conf}.  Second, the pairing
\[ \mQ_4(u,v) = \int_X u\,L_4v + \oint_M \left(B_0^3(u)B_3^3(v) + B_1^3(u)B_2^3(v)\right) \]
is symmetric; in particular, this ensures that the higher order terms in the Paneitz operator and the boundary operators $B_2^3$ and $B_3^3$ cancel after integration by parts.  This claim is the content of Theorem~\ref{thm:boundary}.

To prove the first statement of Theorem~\ref{thm:boundary}, we identify boundary operators of all integer orders between zero and three which are conformally covariant and act on densities of arbitrary weight (cf.\ \cite{GoverPeterson2015}).  After making an additional choice to ensure that $\mQ_4$ is symmetric (cf.\ Remark~\ref{rk:uniqueness}), these operators recover the boundary operators $B_k^3$ in the case of densities of weight $-\frac{n-3}{2}$.  Other choices of weight are relevant to the study of boundary operators associated to other GJMS operators.

The simplest cases are the zeroth- and first-order conformally covariant boundary operators.

\begin{prop}
 \label{prop:oi01}
 Let $(X^{n+1},g)$ be a Riemannian manifold with boundary $M^n=\partial X$ and let $w\in\bR$.  Define $B_{0,w},B_{1,w}\colon C^\infty(X)\to C^\infty(M)$ by
 \begin{align*}
  B_{0,w}u & = u, \\
  B_{1,w}u & = \eta u - \frac{w}{n}Hu .
 \end{align*}
 Given $\sigma\in C^\infty(X)$, set $\hg=e^{2\sigma}g$.  Then
 \begin{align}
  \label{eqn:oi0} \hB_{0,w}(u) & = e^{w\sigma}B_{0,w}\left(e^{-w\sigma}u\right) , \\
  \label{eqn:oi1} \hB_{1,w}(u) & = e^{(w-1)\sigma}B_{1,w}\left(e^{-w\sigma}u\right)
 \end{align}
 for all $u\in C^\infty(M)$.
\end{prop}

\begin{proof}
 It is clear that~\eqref{eqn:oi0} holds.  Observe that the outward-pointing normals $\eta$ and $\heta$ along $M$ defined with respect to $g$ and $\hg$, respectively, are related by $\heta=e^{-\sigma}\eta$.  It follows that
 \begin{equation}
  \label{eqn:Hprime} \hH = e^{-\sigma}\left(H + n\eta\sigma\right),
 \end{equation}
 from which the conclusion immediately follows.
\end{proof}

In order to more easily identify the conformally covariant boundary operators of higher order, we use an observation of Branson~\cite{Branson1985}: Suppose that $D\colon C^\infty(X)\to C^\infty(M)$ is a natural operator which is homogeneous of degree $k$; i.e.\ $D=D_g$ is defined by a formula which is a complete contraction of a tensor polynomial in the Levi-Civita connection, the Riemann curvature tensor, and the outward-pointing normal along $M$, and moreover, $D_{c^2g}=c^{-k}D_g$ for all constants $c>0$.  Given a fixed constant $w\in\bR$ and a fixed function $\sigma\in C^\infty(X)$, define the operator $D^\prime$ by
\[ D^\prime(u) := \left.\frac{d}{dt}\right|_{t=0}\left[ e^{-(w-k)t\sigma}D_{e^{2t\sigma}g}\left(e^{wt\sigma}u\right) \right] . \]
We equivalently write $(Du)^\prime$ for $D^\prime(u)$.  It is straightforward to check that $D^\prime\equiv0$ if and only if
\[ D_{e^{2\sigma}g}(u) = e^{(w-k)\sigma}D_g\left(e^{-w\sigma}u\right) \]
for all $u\in C^\infty(X)$.  Thus to prove that a natural homogeneous operator is conformally covariant, it suffices to show that its ``linearization'' $D^\prime$ is the zero operator for some weight $w$.

The boundary operator $B_2^3$ is a natural operator which is homogeneous of degree two.  To check that it is conformally covariant, we compute the linearization of most (cf.\ Remark~\ref{rk:uniqueness}) natural operators which are homogeneous of degree two.

\begin{lem}
 \label{lem:pre2}
 Let $(X^{n+1},g)$ be a Riemannian manifold with boundary $M^n=\partial X$.  Given $w\in\bR$ and $\sigma\in C^\infty(X)$, it holds that
 \begin{align*}
  \left(\oDelta u\right)^\prime & = (n+2w-2)\lp\onabla u,\onabla\sigma\rp + wu\oDelta\sigma , \\
  \left(\nabla^2u(\eta,\eta)\right)^\prime & = \lp\onabla u,\onabla\sigma\rp + (2w-1)(\eta u)(\eta\sigma) + wu\nabla^2\sigma(\eta,\eta) , \\
  \left(H\eta u\right)^\prime & = n(\eta u)(\eta\sigma) + wHu\eta\sigma, \\
  \left(H^2u\right)^\prime & = 2nHu\eta\sigma, \\
  \left(uP(\eta,\eta)\right)^\prime & = -u\nabla^2\sigma(\eta,\eta), \\
  \left(u\oJ\right)^\prime & = -u\oDelta\sigma
 \end{align*}
 for all $u\in C^\infty(X)$.
\end{lem}

\begin{proof}
 It is well-known that if $\hg=e^{2\sigma}g$, then
 \begin{align*}
  \hP & = P - \nabla^2\sigma + d\sigma\otimes d\sigma - \frac{1}{2}\lv\nabla\sigma\rv^2\,g , \\
  \hnabla^2u & = \nabla^2u - d\sigma\otimes du - du\otimes d\sigma + \lp\nabla\sigma,\nabla u\rp\,g
 \end{align*}
 Combining this, the analogous statements for the conformally rescaled boundary metric $\hh$, and~\eqref{eqn:Hprime} yields the result.
\end{proof}

In particular, we obtain the following family of second-order conformally covariant boundary operators.

\begin{prop}
 \label{prop:oi2}
 Let $(X^{n+1},g)$ be a Riemannian manifold with boundary $M^n=\partial X$ and let $w\in\bR$.  Define $B_{2,w}\colon C^\infty(X)\to C^\infty(M)$ by
 \begin{multline*}
  B_{2,w}u := -\oDelta u + (n+2w-2)\nabla^2u(\eta,\eta) - \frac{(n+2w-2)(2w-1)}{n}H\eta u \\ - w\left(\oJ - (n+2w-2)P(\eta,\eta) - \frac{(n+2w-2)(2w-1)}{2n^2}H^2\right)u .
 \end{multline*}
 Given $\sigma\in C^\infty(X)$, set $\hg=e^{2\sigma}g$.  Then
 \[ \hB_{2,w}(u) = e^{(w-2)\sigma}B_{2,w}\left(e^{-w\sigma}u\right) \]
 for all $u\in C^\infty(X)$.
\end{prop}

\begin{proof}
 Lemma~\ref{lem:pre2} implies that $B_{2,w}^\prime=0$.
\end{proof}

The boundary operator $B_3^3$ is a natural operator which is homogeneous of degree three.  To check that it is conformally covariant, we compute the linearization of most (cf.\ Remark~\ref{rk:uniqueness}) natural operators which are homogeneous of degree three.  For the sake of readability, we separate this computation into four lemmas according to the order of the operator as a differential operator.

\begin{lem}
 \label{lem:pre33}
 Let $(X^{n+1},g)$ be a Riemannian manifold with boundary $M^n=\partial X$.  Given $w\in\bR$ and $\sigma\in C^\infty(X)$, it holds that
 \begin{align*}
  \left(\eta\Delta u\right)^\prime & = (w-2)(\eta\sigma)\oDelta u + (n+3w-3)(\eta\sigma)\nabla^2u(\eta,\eta) \\
   & \quad + (n+2w-1)\lp\onabla\eta u,\onabla\sigma\rp + 2(w-1)H(\eta u)(\eta\sigma) \\
   & \quad + (n+2w-1)\lp\onabla u,\onabla\eta\sigma\rp - 2(n+2w-1)A(\onabla u,\onabla\sigma\rp \\
   & \quad + (n+3w-1)(\eta u)\nabla^2\sigma(\eta,\eta) + w(\eta u)\oDelta\sigma + wu\eta\Delta\sigma , \\
  \left(\oDelta\eta u\right)^\prime & = w(\eta\sigma)\oDelta u + (n+2w-4)\lp\onabla\eta u,\onabla\sigma\rp \\
   & \quad + 2w\lp\onabla u,\onabla\eta\sigma\rp + (w-1)(\eta u)\oDelta\sigma + wu\oDelta\eta\sigma
 \end{align*}
 for all $u\in C^\infty(X)$.
\end{lem}

\begin{proof}
 Lemma~\ref{lem:pre2} yields the formula for $(\oDelta\eta)^\prime$ and also that
 \[ \left(\eta\Delta u\right)^\prime = \eta\left((n+2w-1)\lp\nabla u,\nabla\sigma\rp + wu\Delta\sigma\right) + (w-2)(\eta\sigma)\Delta u . \]
 The result then follows from Lemma~\ref{lem:2derivatives} and Lemma~\ref{lem:derivatives_eta}.
\end{proof}

\begin{lem}
 \label{lem:pre32}
 Let $(X^{n+1},g)$ be a Riemannian manifold with boundary $M^n=\partial X$.  Given $w\in\bR$ and $\sigma\in C^\infty(X)$, it holds that
 \begin{align*}
  \left(\lp A_0,\onabla^2u\rp\right)^\prime & = 2(w-1)A_0(\onabla u,\onabla\sigma) + wu\lp A_0,\onabla^2\sigma\rp, \\
  \left(H\oDelta u\right)^\prime & = n(\eta\sigma)\oDelta u + (n+2w-2)H\lp\onabla u,\onabla\sigma\rp + wHu\oDelta\sigma, \\
  \left(H\nabla^2u(\eta,\eta)\right)^\prime & = n(\eta\sigma)\nabla^2u(\eta,\eta) + H\lp\onabla u,\onabla\sigma\rp \\
   & \quad + (2w-1)H(\eta u)(\eta\sigma) + wHu\nabla^2\sigma(\eta,\eta)
 \end{align*}
 for all $u\in C^\infty(X)$.
\end{lem}

\begin{proof}
 Recall that $A_0$ is conformally covariant: If $\hg=e^{2\sigma}g$, then $\hA_0=e^{\sigma}A_0$.  This observation, \eqref{eqn:Hprime} and Lemma~\ref{lem:pre2} together yield the result.
\end{proof}

\begin{lem}
 \label{lem:pre31}
 Let $(X^{n+1},g)$ be a Riemannian manifold with boundary $M^n=\partial X$.  Given $w\in\bR$ and $\sigma\in C^\infty(X)$, it holds that
 \begin{align*}
  \left(\lp\onabla H,\onabla u\rp\right)^\prime & = n\lp\onabla u,\onabla\eta\sigma\rp - H\lp\onabla u,\onabla\sigma\rp + wu\lp\onabla H,\onabla\sigma\rp, \\
  \left(\oJ\eta u\right)^\prime & = -(\eta u)\oDelta\sigma + w\oJ u\eta\sigma, \\
  \left(P(\eta,\eta)\eta u\right)^\prime & = -(\eta u)\nabla^2\sigma(\eta,\eta) + wuP(\eta,\eta)\eta\sigma , \\
  \left(H^2\eta u\right)^\prime & = 2nH(\eta u)\eta\sigma + wH^2u\eta\sigma
 \end{align*}
 for all $u\in C^\infty(X)$.
\end{lem}

\begin{proof}
 This follows immediately from~\eqref{eqn:Hprime} and Lemma~\ref{lem:pre2}.
\end{proof}

\begin{lem}
 \label{lem:pre30}
 Let $(X^{n+1},g)$ be a Riemannian manifold with boundary $M^n=\partial X$.  Given $w\in\bR$ and $\sigma\in C^\infty(X)$, it holds that
 \begin{align*}
  \left(u\eta J\right)^\prime & = -u\eta\Delta\sigma - 2\left(\oJ+P(\eta,\eta) + \frac{1}{2(n-1)}\lv A_0\rv^2 - \frac{1}{2n}H^2\right)u\eta\sigma , \\
  \left(u\oDelta H\right)^\prime & = nu\oDelta\eta\sigma - Hu\oDelta\sigma + (n-4)u\lp\onabla H,\onabla\sigma\rp, \\
  \left(uH\oJ\right)^\prime & = n\oJ u\eta\sigma - Hu\oDelta\sigma, \\
  \left(uHP(\eta,\eta)\right)^\prime & = nuP(\eta,\eta)\eta\sigma - Hu\nabla^2\sigma(\eta,\eta), \\
  \left(u\lp A_0,\oP\rp\right)^\prime & = -u\lp A_0,\onabla^2\sigma\rp, \\
  \left(uH^3\right)^\prime & = 3nH^2u\eta\sigma , \\
  \left(uH\lv A_0\rv^2\right) & = n\lv A_0\rv^2u\eta\sigma
 \end{align*}
 for all $u\in C^\infty(X)$.
\end{lem}

\begin{proof}
 This follows immediately from Lemma~\ref{lem:gauss_codazzi}, \eqref{eqn:Hprime}, Lemma~\ref{lem:pre2} and the conformal covariance of $A_0$.
\end{proof}

In particular, we obtain the following family of third-order conformally covariant boundary operators.

\begin{prop}
 \label{prop:oi3}
 Let $(X^{n+1},g)$ be a Riemannian manifold with boundary $M^n=\partial X$ and let $w\in\bR$.  Define $B_{3,w}\colon C^\infty(X)\to C^\infty(M)$ by
 \begin{align*}
  B_{3,w}u & := -\eta\Delta u + \frac{n+2w-1}{n+2w-4}\oDelta\eta u + \frac{n+3w-3}{n}H\nabla^2u(\eta,\eta) \\
   & \quad - \frac{n+2w-1}{w-1}\lp A_0,\onabla^2u\rp - \frac{2n+7w-8}{n(n+2w-4)}H\oDelta u \\
   & \quad + \frac{(n-4)(n+2w-1)}{n(n+2w-4)}\lp\onabla H,\onabla u\rp + S_{2,w}\eta u - wT_{3,w}u,
 \end{align*}
 where
 \begin{align*}
  S_{2,w} & = \left(\frac{w-1}{n} - \frac{(n+3w-3)(2w-1)}{2n^2}\right)H^2 - (n+3w-1)P(\eta,\eta) \\
   & \quad - \frac{n-w-1}{n+2w-4}\oJ , \\
  T_{3,w} & = \eta J + \frac{n+2w-1}{n(n+2w-4)}\oDelta H + \frac{n+2w-1}{w-1}\lp A_0,\oP\rp \\
   & \quad - \frac{n+3w-3}{n}HP(\eta,\eta) + \frac{n+5w-7}{n(n+2w-4)}H\oJ \\
   & \quad + \frac{1}{n(n-1)}H\lv A_0\rv^2 + \left(\frac{w-2}{3n^2}-\frac{(n+3w-3)(2w-1)}{6n^3}\right)H^3 .
 \end{align*}
 Given $\sigma\in C^\infty(X)$, set $\hg=e^{2\sigma}g$.  Then
 \[ \hB_{3,w}(u) = e^{(w-3)\sigma}B_{3,w}\left(e^{-w\sigma}u\right) \]
 for all $u\in C^\infty(X)$.
\end{prop}

\begin{proof}
 Lemma~\ref{lem:pre33}, Lemma~\ref{lem:pre32}, Lemma~\ref{lem:pre31}, and Lemma~\ref{lem:pre30} together imply that $B_{3,w}^\prime=0$.
\end{proof}

We are now prepared to prove that the operators $B_k^3$ for $k\in\{0,1,2,3\}$ are boundary operators associated to the Paneitz operator.

\begin{proof}[Proof of Theorem~\ref{thm:boundary}]
 Observe that $B_j^3=B_{j,-\frac{n-3}{2}}$ for $j\in\{0,1,2\}$ and that
 \[ B_3^3 = B_{3-\frac{n-3}{2}} + \frac{1}{2}\lv A_0\rv^2B_{1,-\frac{n-3}{2}} . \]
 It follows from Proposition~\ref{prop:oi01}, Proposition~\ref{prop:oi2} and Proposition~\ref{prop:oi3} that if $\hg=e^{2\sigma}g$, then
 \[ \hB_k^3(u) = e^{-\frac{n-3+2k}{2}\sigma}B_k^3\left(e^{\frac{n-3}{2}\sigma}u\right) \]
 for all $u\in C^\infty(X)$.

 Next, the divergence theorem implies that
 \begin{multline}
  \label{eqn:paneitz_ibp}
  \int_X \left( u\,L_4v - (\Delta u)(\Delta v) + (4P-(n-1)Jg)(\nabla u,\nabla v) - \frac{n-3}{2}Quv\right) \\ = \oint_M \left(u\eta\Delta v - (\Delta v)\eta u + 4uP(\eta,\nabla v) - (n-1)Ju\eta v\right) .
 \end{multline}
 for all $u,v\in C^\infty(X)$.  On the other hand, Lemma~\ref{lem:gauss_codazzi} and Lemma~\ref{lem:2derivatives} imply that
 \begin{align}
  \label{eqn:pibp1} \oint_M Ju\eta v & = \oint_M \left(\oJ + P(\eta,\eta) - \frac{1}{2n}H^2 + \frac{1}{2(n-1)}\lv A_0\rv^2\right)u\eta v, \\
  \label{eqn:pibp2} \oint_M uP(\eta,\nabla v) & = \oint_M \left(uP(\eta,\eta)\eta v + \frac{u}{n-1}\lp\odelta A - \onabla H,\onabla v\rp \right), \\
  \label{eqn:pibp3} \oint_M (\Delta v)\eta u & = \oint_M \left( (\eta u)\nabla^2v(\eta,\eta) - \lp\onabla v,\onabla\eta u\rp + H(\eta u)(\eta v)\right)
 \end{align}
 for all $u,v\in C^\infty(X)$.  Inserting these formulae into~\eqref{eqn:paneitz_ibp} and using the definitions of the boundary operators $B_k^3$ for $k\in\{0,1,2,3\}$ yields~\eqref{eqn:l4_symmetric}.  It is then clear that $\mQ_4$ is a symmetric bilinear form.
\end{proof}

\begin{remark}
 \label{rk:uniqueness}
 From invariant theory, it is easy to check that the families $B_{0,w}$ and $B_{1,w}$ are the unique (up to constant multiples) families of conformally covariant operators which are natural and homogeneous of degree zero and one, respectively.  However, the families $B_{2,w}$ and $B_{3,w}$ are not unique.  The scalar $\lv A_0\rv^2$ is conformally covariant and homogeneous of degree two.  By invariant theory and Lemma~\ref{lem:2derivatives}, the span of $B_{2,w}$ and $\lv A_0\rv^2$ gives the set of all natural conformally covariant boundary operators which are homogeneous of degree two.  Likewise, the operator $\lv A_0\rv^2B_{1,w}$ and the scalars $\tr A_0^3$ and $\lp W(\eta,\cdot,\eta,\cdot),A_0\rp$ are conformally covariant and homogeneous of degree three.  By invariant theory and Lemma~\ref{lem:3derivatives}, the span of $B_3^3$ and these three operators gives the set of all natural conformally covariant boundary operators which are homogeneous of degree three.

 The requirement that the operators $B_k^3$ for $k\in\{0,1,2,3\}$ be boundary operators for the Paneitz operator imposes, in addition to conformal covariance, the requirement that $\mQ_4$ be symmetric.  This requirement imposes three additional constraints.  In order to absorb the term
 \[ \int_X \left( u\Delta^2v - (\Delta u)(\Delta v)\right) = \oint_M \left( u\eta\Delta v - (\oDelta v + \nabla^2v(\eta,\eta) + H\eta v)\eta u\right) \]
 arising in integration by parts, it must hold that
 \begin{align*}
  B_0^3(u)B_3^3(v) & = -u\eta\Delta v + \text{(lower order terms)}, \\
  B_1^3(v)B_2^3(v) & = (\eta u)\nabla^2v(\eta,\eta) + \text{(lower order terms)},
 \end{align*}
 where ``lower order terms'' means terms with a total number of normal derivatives less than three.  To absorb the term involving $\oint\lv A_0\rv^2u\eta v$ in~\eqref{eqn:pibp1} requires that there is a $k\in\bR$ such that
 \[ B_0^3(u)B_3^3(v) + B_1^3(u)B_2^3(v) = \frac{1}{2}\lv A_0\rv^2\left((1+k)u\eta v + kv\eta u\right) \]
 up to terms not involving $\lv A_0\rv^2\eta$.

 In summary, if we normalize $B_0^3$ and $B_1^3$ to have leading normal order terms $1$ and $\eta$, respectively, then there are three remaining freedoms in how we define the boundary operators: We can add a linear combination of $\tr A_0^3$ and $\lp W(\eta,\cdot,\eta,\cdot),A_0\rp$ to $B_3^3$; and we can add a constant multiple of $\lv A_0\rv^2$ to $B_2^3$ while adding the same constant multiple of $\lv A_0\rv^2B_1^3$ to $B_3^3$.
\end{remark}

One simple consequence of Theorem~\ref{thm:boundary} is the geometric interpretation of positive smooth minimizers of the energy $\mE_4$.

\begin{proof}[Proof of Proposition~\ref{prop:l4_geometric_euler}]
 Consider first the case that $u\in C^\infty(X)$ is a positive minimizer of $Y_{4,2}(X,M)$.  Therefore $u$ is a critical point of $\mE_4$ when restricted to the class of smooth functions $w\in C^\infty(X)$ with $\oint \lv B_0^3w\rv^{\frac{2n}{n-3}}=1$.  Let $u_t$ be a one-parameter family of such functions which is $C^1$ in $t\in(-\varepsilon,\varepsilon)$ and satisfies $u_0=u$.  Using Theorem~\ref{thm:boundary}, we compute that
 \begin{align*}
  0 & = \frac{1}{2}\left.\frac{d}{dt}\right|_{t=0}\mE_4(u_t) = \int_X v\,L_4u + \oint_M \left[ B_0^3(v)B_3^3(u) + B_1^3(v)B_2^3(u) \right] , \\
  0 & = \frac{1}{2}\left.\frac{d}{dt}\right|_{t=0}\int_M u_t^{\frac{2n}{n-3}} = \frac{n}{n-3}\oint_M vu^{\frac{n+3}{n-3}}
 \end{align*}
 for $v=\frac{\partial u_t}{\partial t}\rv_{t=0}$. Thus there is a constant $\lambda\in\bR$ such that
 \begin{equation}
  \label{eqn:l4_geometric_euler}
  \begin{cases}
   L_4u = 0, \\
   B_2^3u = 0, \\
   B_3^3u = \lambda u^{\frac{n+3}{n-3}} .
  \end{cases}
 \end{equation}
 From the assumption $\mE_4(u)=Y_{4,2}(X,M)$ we compute that $\lambda=Y_{4,2}(X,M)$.  The final conclusion for $\hg=u^{\frac{4n}{n-3}}g$ follows from Theorem~\ref{thm:boundary}.

 Consider next the case that $u\in C^\infty(X)$ is a positive minimizer of $Y_{4,1}(X,M)$.  Therefore $u$ is a critical point of $\mE_4$ when restricted to the class of smooth functions $w\in C^\infty(X)$ with $B_1^3w=0$ and $\oint\lv B_0^3w\rv^{\frac{2n}{n-3}}=1$.  Arguing as above, we deduce that
 \[ \begin{cases}
     L_4u=0, \\
     B_1^3u=0, \\
     B_3^3u = Y_{4,1}(X,M) u^{\frac{n+3}{n-3}} .
    \end{cases} \]
 The final conclusion for $\hg=u^{\frac{4n}{n-3}}g$ again follows from Theorem~\ref{thm:boundary}.
\end{proof}

\begin{remark}
 The Euler equation~\eqref{eqn:l4_geometric_euler} holds for any boundary volume-constrained critical point $u$ of $\mE_4$.  By changing the constraints, it is clear that the energy functional $\mE_4$ gives a variational approach to simultaneously prescribing the interior $Q$-curvature, the $T_3^3$-curvature of the boundary, and one of the mean curvature and the $T_2^3$-curvature of the boundary.
\end{remark}
\section{Boundary operators and pseudodifferential operators}
\label{sec:asymptotics}

From their explicit expressions, it is clear that the Paneitz operator $L_4$ is properly elliptic and that the boundary operators $B_k^3$ are of the form $B_k^3=(\nabla_\eta)^k + \text{(lower order terms)}$ for $k\in\{0,1,2,3\}$, where ``lower order terms'' means terms which involve fewer than $k$ derivatives in the normal direction.  It follows that the boundary value problems $(L_4;(B_0^3,B_1^3))$, $(L_4;(B_0^3,B_2^3))$ and $(L_4;(B_1^3,B_3^3))$ all satisfy the Lopatinskii--Shapiro conditions (cf.\ \cite{AgmonDouglisNirenberg1964,BransonGover2001}).  Moreover, the symmetry of the bilinear form $\mQ_4$ implies that each of these boundary value problems is formally self-adjoint.  Therefore, so long as the null space of the given boundary value problem is trivial, to any pair of smooth initial data there exists a unique smooth solution $u\in\ker L_4$.  For example, if $\ker L_4\cap \ker B_0^3 \cap \ker B_1^3 = \{0\}$, then given $f,\psi\in C^\infty(M)$, there exists a unique function $u_{f,\psi}\in C^\infty(X)$ such that
\begin{equation}
 \label{eqn:unique_01}
 \begin{cases}
  L_4u_{f,\psi} = 0, \\
  B_0^3u_{f,\psi} = f, \\
  B_1^3u_{f,\psi} = \psi .
 \end{cases}
\end{equation}

The solution to~\eqref{eqn:unique_01} can be constructed variationally.  This is simplest when $\lambda_1(L_4)>0$, in which case the solution to~\eqref{eqn:unique_01} minimizes the energy $\mE_4$ within the class of functions $u\in C^\infty(X)$ with fixed trace $(f,\psi)=(B_0^3u,B_1^3u)$.

\begin{proof}[Proof of Proposition~\ref{prop:extension}]
 Let
 \[ \mC_{f,\psi} = \left\{ u \in C^\infty(X) \suchthat B_0^3u=f, B_1^3u=\psi \right\} . \]
 It follows from Theorem~\ref{thm:boundary} that critical points of $\mE_4\colon\mC_{f,\psi}\to\bR$ satisfy~\eqref{eqn:unique_01}.

 Fix $u\in\mC_{f,\psi}$.  Then $\mC_{f,\psi}=u+\mC_{0,0}$.  Given $v\in\mC_{0,0}$, Theorem~\ref{thm:boundary} implies that
 \[ \mE_4(u+v) = \mE_4(v) + 2\int_X v\,L_4u + \mE_4(u) . \]
 Using the assumption $\lambda_1(L_4)>0$, we find that
 \[ \mE_4(u+v) \geq \mE_4(u) - \frac{1}{\lambda_1(L_4)}\int_X\left(L_4u\right)^2 . \]
 Therefore $\mE_4$ is bounded below on $\mC_{f,\psi}$.  By taking a minimizing sequence, we thus obtain a weak solution $u$ of~\eqref{eqn:unique_01}.  Since $(L_4;(B_0^3,B_1^3))$ satisfies the Lopatinskii--Shapiro conditions, it holds that $u\in C^\infty(X)$.

 Finally, if $u_1,u_2$ both satisfy~\eqref{eqn:unique_01}, then $w=u_1-u_2$ is such that $w\in\ker L_4\cap \ker B_0^3\cap\ker B_1^3$, and hence $w\equiv0$.
\end{proof}

Suppose now that $(X^{n+1},g)$ is a compact Riemannian manifold with boundary $M^n=\partial X$ such that the null space of the boundary value problem $(L_4;(B_0^3,B_1^3))$ is trivial.  Hence the operators $\mB_1^3,\mB_3^3\colon C^\infty(M)\to C^\infty(M)$ given by
\begin{align*}
 \mB_1^3\psi & = \frac{1}{2}B_2^3u_{0,\psi}, \\
 \mB_3^3f & = \frac{1}{2}B_3^3u_{f,0}
\end{align*}
are well-defined.  In fact, $\mB_1^3$ and $\mB_3^3$ are formally self-adjoint conformally covariant first- and third-order pseudodifferential operators, respectively:

\begin{prop}
 \label{prop:mB_conf}
 Let $(X^{n+1},g)$ be a compact Riemannian manifold with boundary $M^n=\partial X$ such that the null space of $(L_4;(B_0^3,B_1^3))$ is trivial.  Then $\mB_1^3$ and $\mB_3^3$ are formally self-adjoint conformally covariant pseudodifferential operators with principal symbols $\sigma(\mB_k^3)=\sigma((-\Delta)^{k/2})$ for $k\in\{1,3\}$.  Indeed, given $w\in C^\infty(M)$, set $\hh=e^{2w}h$.  Then
 \begin{equation}
  \label{eqn:psido_conf_cov}
  \hmB_k^3(\psi) = e^{-\frac{n+k}{2}w}\mB_k^3\left(e^{\frac{n-k}{2}w}\psi\right)
 \end{equation}
 for $k\in\{1,3\}$ and all $f,\psi\in C^\infty(M)$.
\end{prop}

\begin{proof}
 Let $f_1,f_2\in C^\infty(M)$.  For $j\in\{1,2\}$, let $u_j=u_{f_j,0}$ be as in~\eqref{eqn:unique_01}.  It follows from Theorem~\ref{thm:boundary} that
 \[ \oint_M f_1\mB_3^3f_2 = \frac{1}{2}\mE_4(u_1,u_2) = \oint_M f_2\mB_3^3 f_1 . \]
 That $\mB_1^3$ is formally self-adjoint follows similarly.

 Let $v\in C^\infty(X)$ be such that $B_0^3v=w$.  Set $\hg=e^{2v}g$.  Given $f,\psi\in C^\infty(M)$, let $u_{f,\psi},\hu_{\hf,\hpsi}\in C^\infty(X)$ be such that
 \begin{align*}
  L_4u_{f,\psi}&=0, & B_0^3u_{f,\psi}&=f, & B_1^3u_{f,\psi}&=\psi, \\
  \hL_4\hu_{\hf,\psi}&=0, & \hB_0^3\hu_{\hf,\hpsi}&=e^{-\frac{n-3}{2}w}f, & \hB_1^3\hu_{\hf,\hpsi}&=e^{-\frac{n-1}{2}w}\psi .
 \end{align*}
 It follows from Theorem~\ref{thm:boundary} that
 \[ \hu_{\hf,\hpsi} = e^{-\frac{n-3}{2}v}u_{f,\psi} . \]
 Applying Theorem~\ref{thm:boundary} in the special cases $f=0$ and $\psi=0$ yields the conformal covariance of $\mB_1^3$ and $\mB_3^3$, respectively.

 The claim that $\mB_k^3$ are pseudodifferential operators for $k\in\{1,3\}$ follows as in~\cite{BransonGover2001}: From the homogeneity of $\mB_k^3$, which is an immediate consequence of~\eqref{eqn:psido_conf_cov}, and the fact that the Paneitz operator and its boundary operators are all $O(TM)$-invariant, it follows that $\sigma(\mB_k^3)(x,\xi)=c_k(x)\lv\xi\rv^k$ for some $c_k\in C^\infty(M)$.  Since the principal symbol of the Paneitz operator is that of the bi-Laplacian $\Delta^2$, the principal symbols of the extension operators $f\mapsto u_{f,0}$ and $\psi\mapsto u_{0,\psi}$ depend only on $\lv\xi\rv$, and hence $c_k$ is a universal constant~\cite{Kumano-go1981}.  It follows from Theorem~\ref{thm:asymptotics} and the identification of the principal symbol of the fractional GJMS operators~\cite{GrahamZworski2003} that $c_k=1$.
\end{proof}

\begin{remark}
 Note that our proof of conformal covariance shows that the operators $\mB_k^3$ for $k\in\{1,3\}$ depend only on $[g]$ and a choice of metric $h$ in the induced conformal class on $M$.
\end{remark}

There is another way one can define first- and third-order conformally covariant pseudodifferential operators on the boundary of a compact Riemannian manifold $(X^{n+1},g)$, namely via scattering theory~\cite{GrahamZworski2003}.  Let $g_+\in[g]$ be an asymptotically hyperbolic metric in the interior $X_0$ of $X$ and let $r\geq0$ be the geodesic defining function associated to $h$.  Suppose that $\frac{n^2-1}{4},\frac{n^2-9}{4}\not\in\sigma_{pp}(-\Delta_{g_+})$ and that
\begin{equation}
 \label{eqn:asympt_h}
 \left.\frac{\partial}{\partial r}\right|_{r=0}\log\det(h_r^{-1}h_r) = 0 = \left.\frac{\partial^3}{\partial r^3}\right|_{r=0}\log\det(h_r^{-1}h_r) ,
\end{equation}
where $g_+=r^{-2}(dr^2+h_r)$ near $M$ for $h_r$ a one-parameter family of Riemannian metrics on $M$.  Fix $f,\psi\in C^\infty(M)$.  By the spectral assumption, there are unique functions $v_f,v_\psi\in C^\infty(X_0)$ such that
\begin{equation}
 \label{eqn:poisson}
 \begin{split}
  -\Delta v_f - \frac{n^2-9}{4}v_f = 0,& \qquad -\Delta v_\psi - \frac{n^2-1}{4}v_\psi = 0, \\
  \lim_{r\to 0} r^{-\frac{n-3}{2}}v_f = f, & \qquad \lim_{r\to 0} r^{-\frac{n-1}{2}}v_\psi = \psi .
 \end{split}
\end{equation}
From the condition~\eqref{eqn:asympt_h}, we compute that there are functions $G_f,G_\psi\in C^\infty(M)$ such that
\begin{align}
 \label{eqn:asympt_P1} v_\psi & = r^{\frac{n-1}{2}}\left(\psi + G_\psi r - \frac{r^2}{2}\left(\oDelta\psi + \frac{n-1}{4}(h^{ij}h_{ij}^\prime)^\prime\psi \right) + O(r^3) \right) , \\
 \label{eqn:asympt_P3} v_f & = r^{\frac{n-3}{2}}\left(f + \frac{r^2}{2}\left(\oDelta f + \frac{n-3}{4}(h^{ij}h_{ij}^\prime)^\prime f\right) + G_f r^3 + O(r^4) \right) ,
\end{align}
In other words, the condition~\eqref{eqn:asympt_h} implies that the scattering operators $S\left(\frac{n+1}{2}\right)$ and $S\left(\frac{n+3}{2}\right)$ have no residues (cf.\ \cite{GrahamZworski2003,GuillarmouGuillope2007}).  Hence the operators $P_1,P_3\colon C^\infty(M)\to C^\infty(M)$ given by
\begin{align*}
 P_1\psi & := -G_\psi, & P_3f & := 3G_f ,
\end{align*}
are well-defined; indeed, $P_1$ and $P_3$ are first- and third-order conformally covariant pseudodifferential operators with principal symbol that of $(-\Delta)^{1/2}$ and $(-\Delta)^{3/2}$, respectively~\cite{GrahamZworski2003}.

It is natural to ask if $\mB_1=P_1$ and $\mB_3=P_3$ in the setting of the previous paragraph.  While we cannot answer this question in general, we can answer it in the case when $(X^{n+1},g)$ is a compactification of a Poincar\'e--Einstein manifold.  In fact, we prove a stronger result in this setting: the fractional GJMS operators $P_1$ and $P_3$ can be computed directly from the boundary operators $B_k^3$ when acting on Paneitz-harmonic functions.

\begin{proof}[Proof of Theorem~\ref{thm:asymptotics}]
 By Theorem~\ref{thm:asymptotics} and the conformal covariance of the fractional GJMS operators $P_1$ and $P_3$, it suffices to prove the result for the metric
 \begin{align*}
  g & = r^2g_+ = dr^2 + h_r, \\
  h_r & = h - r^2\oP + O(r^3) ,
 \end{align*}
 where $r$ is a geodesic defining function.  Moreover, the coefficient of $r^3$ in the expansion of $h_r$ is tracefree; indeed, it vanishes for dimensions at least four~\cite{FeffermanGraham2012}.  In particular, $A$, $P(\eta,\eta)$ and $\eta J$ all vanish identically along $M$.
 
 Let $u\in C^\infty(X)$ satisfy $L_4u=0$.  Set $f=B_0^3u$ and $\psi=B_1^3u$.  Let $v_f$ and $v_\psi$ be as in~\eqref{eqn:poisson}.  It follows from~\eqref{eqn:asympt_P1} and~\eqref{eqn:asympt_P3} that
 \begin{align*}
  r^{\frac{n-3}{2}}v_f & = f - \frac{r^2}{2}\left(-\oDelta f + \frac{n-3}{2}\oJ f\right) + \frac{r^3}{3}P_3f + O(r^4) , \\
  r^{\frac{n-3}{2}}v_\psi & = r\psi - r^2P_1\psi + \frac{r^3}{2}\left(-\oDelta\psi + \frac{n-1}{2}\oJ\psi \right) + O(r^4) ,
 \end{align*}
 while the uniqueness of solutions to the problem $(L_4;(B_0^3,B_1^3))$, the factorization~\eqref{eqn:paneitz_einstein}, and conformal covariance implies that $u=\left(v_f-v_\psi\right)r^{\frac{n-3}{2}}$.  Therefore
 \begin{multline*}
  u = f - \psi r + \left(P_1\psi - \frac{1}{2}\left(-\oDelta f + \frac{n-3}{2}\oJ f\right)\right)r^2 \\ + \left(\frac{1}{3}P_3f  - \frac{1}{2}\left(-\oDelta\psi + \frac{n-1}{2}\oJ\psi\right)\right)r^3 + O(r^4) .
 \end{multline*}
 In particular,
 \begin{align*}
  \nabla^2u(\eta,\eta) & = 2P_1\psi + \oDelta f - \frac{n-3}{2}\oJ f, \\
  -\eta\Delta u & = 2P_3 + 2\oDelta\psi - \frac{3n-5}{2}\oJ\psi .
 \end{align*}
 Inserting these identities into the definitions of $B_2^3$ and $B_3^3$ yields
 \[ B_2^3 = 2P_1\psi, \qquad B_3^3u = 2P_3f . \qedhere \]
\end{proof}

An immediate corollary of Theorem~\ref{thm:asymptotics} is that for compactifications of Poincar\'e--Einstein manifolds, if $u\in\ker L_4$, then $B_2^3u=2\mB_1^3B_1^3u$ and $B_3^3u=2\mB_3^3B_0^3u$.  This decoupling is helpful when studying Paneitz-harmonic functions, and so it is useful to characterize more generally when it occurs.

\begin{lem}
 \label{lem:mB_uniqueness}
 Let $(X^{n+1},g)$ be a compact Riemannian manifold with boundary $M^n=\partial X$ and suppose that $\ker L_4\cap\ker B_0^3\cap\ker B_1^3=\{0\}$.
 \begin{enumerate}
  \item $B_2^3u_{f,0}=0$ for all $f\in C^\infty(M)$ if and only if $B_2^3u_{f,\psi}=2\mB_1^3\psi$ for all $f,\psi\in C^\infty(M)$.
  \item $B_3^3u_{0,\psi}=0$ for all $\psi\in C^\infty(M)$ if and only if $B_3^3u_{f,\psi}=2\mB_3^3f$ for all $f,\psi\in C^\infty(M)$.
 \end{enumerate}
\end{lem}

\begin{proof}
 Since $\ker L_4\cap\ker B_0^3\cap\ker B_1^3=\{0\}$, given $f,\psi\in C^\infty(M)$, there is a unique solution $u_{f,\psi}\in C^\infty(X)$ of the boundary value problem~\eqref{eqn:unique_01}.  Since $L_4$ and the boundary operators $B_k^3$ are all linear, it holds that $u_{f,\psi}=u_{f,0}+u_{0,\psi}$.  In particular, if $B_2^3u_{f,0}=0$ for all $f\in C^\infty(M)$, then
 \[ B_2^3u_{f,\psi} = B_2^3u_{0,\psi} = 2\mB_1^3\psi . \]
 The other conclusions follow similarly.
\end{proof}
\section{Sharp Sobolev trace inequalities}
\label{sec:sobolev}

Proposition~\ref{prop:extension} and Theorem~\ref{thm:asymptotics} can be applied to establish sharp Sobolev trace inequalities involving $W^{2,2}$ on standard models.  We focus here on Euclidean upper half space $(\bR_+^{n+1},dx^2)$ and the round upper hemisphere $(S_+^{n+1},d\theta^2)$, though similar results on other conformally equivalent spaces readily follow from Theorem~\ref{thm:boundary}.

First note that Proposition~\ref{prop:extension} gives a norm computation for the Sobolev trace embedding
\[ W^{2,2}(X^{n+1}) \hookrightarrow H^{3/2}(M) \oplus H^{1/2}(M) \]
and an identification of a right inverse.  To make this completely rigorous, \eqref{eqn:extension_energy} must be rewritten with the $L^2$-norm of $\nabla^2u$ rather than $\Delta u$, as otherwise the estimates degenerate.  For example, if $\psi=(-\Delta)^{1/2}f$ on $\bR^n$, then the harmonic extension $u_f$ of $f$ to $\bR_+^{n+1}$ is also the solution $u_{f,\psi}$.  The correct formulation follows from Reilly's formula~\cite{Reilly1977}, which we recast here in terms of $B_0^3$ and $B_1^3$ as follows.

\begin{lem}
 \label{lem:reilly}
 Let $(X^{n+1},g)$ be a compact Riemannian manifold with boundary $M^n=\partial X$.  Let $u\in C^\infty(X)$ and set $f=B_0^3u$ and $\psi=B_1^3u$.  Then
 \begin{align*}
  \int_X \lv\nabla^2u\rv^2 & = \int_X \left[ \left(\Delta u\right)^2 - (n-1)P(\nabla u,\nabla u) - J\lv\nabla u\rv^2 \right] \\
  & \quad - \oint_M \biggl[ H\psi^2 - 2\lp\onabla\psi,\onabla f\rp - (n-3)H^2f\psi \\
   & \qquad + \left(A+(n-3)Hh\right)(\onabla f,\onabla f) - \frac{n-3}{2}\left(\oDelta H - \frac{n-3}{2}H^3\right)f^2 \biggr] .
 \end{align*}
\end{lem}

\begin{proof}
 Recall Reilly's formula~\cite{Reilly1977}
 \begin{multline*}
  \int_X \left[ \lv\nabla^2u\rv^2 - \left(\Delta u\right)^2 + (n-1)P(\nabla u,\nabla u) + J\lv\nabla u\rv^2 \right] \\ = \oint_M \left[ 2\lp\onabla\eta u,\onabla u\rp - A(\onabla u,\onabla u) - H(\eta u)^2 \right] .
 \end{multline*}
 The conclusion follows from the definitions of $B_0^3$ and $B_1^3$.
\end{proof}

\begin{cor}
 \label{cor:trace}
 On Euclidean upper half space $(\bR_+^{n+1},dx^2)$,
 \[ \int_{\bR_+^{n+1}} \lv\nabla^2u\rv^2 \geq 2\int_{\bR^n} \left[ f(-\Delta)^{3/2}f - \lp\onabla f,\onabla\psi\rp + \psi(-\Delta)^{1/2}\psi \right] \]
 for all $u\in C_0^\infty(\bR_+^{n+1})$, where $f(x)=u(x,0)$ and $\psi(x)=-u_y(x,0)$.  Moreover, equality holds if and only if $\Delta^2u=0$.
\end{cor}

\begin{proof}
 From Lemma~\ref{lem:reilly} we compute that
 \[ \int_{\bR_+^{n+1}}\lv\nabla^2u\rv^2 = \int_{\bR_+^{n+1}} (\Delta u)^2 + 2\oint_{\bR^n}\lp\onabla\psi,\onabla f\rp . \]
 On the other hand, Proposition~\ref{prop:extension} and Theorem~\ref{thm:asymptotics} imply that
 \begin{equation}
  \label{eqn:almost_trace}
  \int_{\bR_+^{n+1}} (\Delta u)^2 \geq \oint_{\bR^n} \left[ 2f(-\Delta)^{3/2}f - 4\lp\onabla f,\onabla\psi\rp + 2\psi(-\Delta)^{1/2}\psi \right]
 \end{equation}
 with equality if and only if $\Delta^2u=0$.  Combining these facts yields the result.
\end{proof}

Our second goal is to prove Theorem~\ref{thm:sobolev}.  This gives a sharp Sobolev trace inequality which gives a norm computation for the embedding
\[ W^{2,2}(\bR_+^{n+1}) \hookrightarrow L^{\frac{2n}{n-3}}(\bR^n) \oplus L^{\frac{2n}{n-1}}(\bR^n) . \]
The embedding can be realized from the trace map and the embeddings $H^\gamma(\bR^n)\hookrightarrow L^{\frac{2n}{n-2\gamma}}(\bR^n)$ for $\gamma\in\{1/2,3/2\}$.  The proof of Theorem~\ref{thm:sobolev} also proceeds in this manner: We apply Proposition~\ref{prop:extension} and Theorem~\ref{thm:asymptotics} to relate the energy $\mE_4$ to the energies of the fractional GJMS operators $P_1$ and $P_3$, and then apply the sharp fractional Sobolev inequalities~\cite{Beckner1993,Lieb1983}.

\begin{proof}[Proof of Theorem~\ref{thm:sobolev}]
 Let $u\in C^\infty(S_+^{n+1})$ and set $f=B_0^3u$ and $\psi=B_1^3u$.  It follows from Proposition~\ref{prop:extension} and Theorem~\ref{thm:asymptotics} that
 \[ \frac{1}{2}\mE_4(u) \geq \oint_{S^n} f\,P_3f + \oint_{S^n} \psi\,P_1\psi . \]
 On the other hand, the fractional Sobolev inequalities (cf.\ \cite{Beckner1993,Lieb1983}) state that
 \[ \oint_{S^n} w\,P_{2\gamma}w \geq C_{2\gamma}\left(\oint_{S^n} \lv w\rv^{\frac{2n}{n-2\gamma}}\right)^{\frac{n-2\gamma}{n}} \]
 for all $w\in C^\infty(S^n)$ and all $\gamma\in(0,n/2)$, where
 \[ C_{2\gamma} = 2^{2\gamma}\pi^\gamma\frac{\Gamma\left(\frac{n+2\gamma}{2}\right)}{\Gamma\left(\frac{n-2\gamma}{2}\right)}\left(\frac{\Gamma(n/2)}{\Gamma(n)}\right)^{\frac{2\gamma}{n}} . \]
 Moreover, equality holds if and only if there are constants $a_1,a_2\in\bR$ and points $\xi_1,\xi_2\in B_1(0)\subset\bR^{n+1}$ such that $f(x)=a_1(1+x\cdot\xi_1)^{-\frac{n-3}{2}}$ and $\psi(x)=a_2(1+\xi\cdot\xi_2)^{-\frac{n-1}{2}}$.  The conclusion now follows from the fact that $S_+^{n+1}$ is Einstein with $P=\frac{1}{2}g$ and the boundary $S^n=\partial S_+^{n+1}$ is totally geodesic.
\end{proof}
\section{Four-manifolds with boundary}
\label{sec:3d}

Consider now the special case of a Riemannian four-manifold $(X^4,g)$ with boundary $M^3=\partial X$.  In this case the boundary operators $B_k^3\colon C^\infty(X)\to C^\infty(M)$ from Theorem~\ref{thm:boundary} are given by
\begin{align*}
 B_0^3u & = u, \\
 B_1^3u & = \eta u, \\
 B_2^3u & = -\oDelta u + \nabla^2u(\eta,\eta) + \frac{1}{3}H\eta u, \\
 B_3^3u & = -\eta\Delta u - 2\oDelta\eta u + 2\lp A_0,\onabla^2u\rp - \frac{2}{3}H\oDelta u + \frac{2}{3}\lp\onabla H, \onabla u \rp + S_2^3\eta u
\end{align*}
for all $u\in C^\infty(X)$, where
\[ S_2^3 = 2\oJ - 2P(\eta,\eta) - \frac{1}{3}H^2 + \frac{1}{2}\lv A_0\rv^2 . \]
For each $k\in\{1,2,3\}$, the kernel of the boundary operator $B_k^3$ constants the constant functions, and thus the ``Branson trick'' \cite{Branson1995} implies that the scalar invariants
\begin{align*}
 T_1^3 & = \frac{1}{3}H, \\
 T_2^3 & = \oJ - P(\eta,\eta) + \frac{1}{18}H^2 , \\
 T_3^3 & = \eta J - \frac{2}{3}\oDelta H - 2\lp A_0,\oP\rp + \frac{4}{3}H\oJ + \frac{1}{3}H\lv A_0\rv^2 - \frac{2}{27}H^3
\end{align*}
act as the $Q$-curvatures associated to the respective operators $B_k^3$.  This can also be verified directly.

\begin{prop}
 \label{prop:l4_conf_crit}
 Let $(X^4,g)$ be a Riemannian four-manifold with boundary $M^3=\partial X$.  Let $\sigma\in C^\infty(X)$ and consider the metric $\hg=e^{2\sigma}g$.  For each $k\in\{1,2,3\}$, it holds that
 \[ e^{k\sigma}\hT_k^3 = T_k^3 + B_k^3\sigma . \]
\end{prop}

\begin{proof}
 \eqref{eqn:Hprime} establishes the claim when $k=1$.  By Theorem~\ref{thm:boundary}, both $B_2^3$ and $B_3^3$ are conformally covariant, and thus it suffices to show that
 \[ \left.\frac{d}{dt}\right|_{t=0}e^{tk\sigma}T_k^3[e^{2t\sigma}] = B_k^3\sigma \]
 for $k\in\{2,3\}$ (cf.\ Section~\ref{sec:operators}).  This follows immediately from Lemma~\ref{lem:pre2} when $k=2$ and from Lemma~\ref{lem:pre30} when $k=3$.
\end{proof}

It follows that $\int Q_4+\oint T_3^3$ is a conformal invariant of a compact Riemannian four-manifold with boundary (cf.\ \cite{ChangQing1997a}).

\begin{cor}
 \label{cor:boundary_q4_invariant}
 Let $(X^4,g)$ be a compact Riemannian manifold with boundary $M^3=\partial X$.  Let $\sigma\in C^\infty(X)$ and consider the metric $\hg=e^{2\sigma}g$.  Then
 \[ \int_X \hQ_4 + \oint_M \hT_3^3 = \int_X Q_4 + \oint_M T_3^3 . \]
\end{cor}

\begin{proof}
 Proposition~\ref{prop:l4_conf_crit} and the transformation formula for the $Q$-curvature imply that
 \[ \int_X\hQ_4\dvol_{\hg} + \oint_M \hT_3^3\dvol_{\hh} = \int_X \left(Q_4 + P_4\sigma\right)\dvol_g + \oint_M \left(T_3^3 + B_3^3\sigma\right)\dvol_h . \]
 Since $P_41=0$ and $B_k^31=0$ for $k\in\{1,2,3\}$, Theorem~\ref{thm:boundary} yields the conclusion.
\end{proof}

It is instructive to compare our boundary operators to the operators defined by Chang and Qing~\cite{ChangQing1997a}.  To that end, define $P_3,P_3^2\colon C^\infty(X)\to C^\infty(M)$ by
\begin{align*}
 P_3(u) & = -\frac{1}{2}\eta\Delta u - \oDelta\eta u + \lp A-\frac{2}{3}H\og,\onabla^2u\rp + \frac{1}{3}\lp\onabla H,\onabla u\rp - \left(\Ric(\eta,\eta)-2J\right)\eta u, \\
 P_3^2(u) & = H\oDelta u - H\nabla^2u(\eta,\eta) - \left(3\Ric(\eta,\eta) - 6J - \frac{1}{3}H^2\right)\eta u
\end{align*}
and set
\[ T = \frac{1}{2}\eta J - \frac{1}{3}\oDelta H + JH - \lp R(\eta,\cdot,\eta,\cdot), A\rp - \frac{1}{3}\tr A^3 + \frac{1}{9}H^3 . \]
These operators are defined exactly as in~\cite{ChangQing1997a}, and differ from our boundary operators by local conformal invariants (cf. Remark~\ref{rk:uniqueness}):

\begin{lem}
 \label{lem:cf_changqing_local}
 Let $(X^4,g)$ be a Riemannian manifold with boundary $M^3=\partial X$.  Then
 \begin{align*}
  P_3 & = \frac{1}{2}B_3^3, \\
  P_3^2 & = -3T_1^3B_2^3 + 3T_2^3B_1^3 + \frac{3}{4}\lv A_0\rv^2B_1^3, \\
  T & = \frac{1}{2}T_3^3 - 2\lp W(\eta,\cdot,\eta,\cdot), A_0\rp - \frac{4}{3}\tr A_0^3 .
 \end{align*}
\end{lem}

\begin{proof}
 Since $M$ is three-dimensional, it holds that
 \begin{align*}
  \lv A\rv^2 & = \lv A_0\rv^2 + \frac{1}{3}H^2 , \\
  \tr A^3 & = \tr A_0^3 + H\lv A_0\rv^2 + \frac{1}{9}H^3 . 
 \end{align*}
 On the other hand, it follows from Lemma~\ref{lem:gauss_codazzi} that
 \begin{align*}
  \lp R(\eta,\cdot,\eta,\cdot), A \rp & = 2\lp W(\eta,\cdot,\eta,\cdot), A_0\rp + \lp\oP,A_0\rp \\
   & \quad + \frac{1}{3}H\oJ + HP(\eta,\eta) - \frac{1}{18}H^3 - \frac{1}{4}H\lv A_0\rv^2 + \tr A_0^3, \\
  \Ric(\eta,\eta) - 2J & = P(\eta,\eta) - \oJ + \frac{1}{6}H^2 - \frac{1}{4}\lv A_0\rv^2 .
 \end{align*}
 Inserting these formulae into the definitions of $P_3,P_3^2,T$ yields the result.
\end{proof}

Proposition~\ref{prop:l4_conf_crit} states that $T_3^3$ can be regarded as the $Q$-curvature associated to $B_3^3$.  By choosing the interior metric appropriately, it also gives rise to a $Q$-curvature associated $\mB_3^3$.

\begin{defn}
 Let $(X^4,g)$ be a compact Riemannian manifold with boundary $M^3=\partial X$ and suppose that $\ker L_4\cap\ker B_0^3\cap\ker B_1^3=\{0\}$.  The $\mT$-curvature of $(M^3,h)$ is defined by
 \[ \mT_3^3 := \frac{1}{2}\hT_3^3, \]
 where $\hT_3^3$ is determined by $\hg=e^{2v}g$ for $v$ the unique solution of
 \begin{equation}
  \label{eqn:mT_ext}
  \begin{cases}
   P_4v + Q_4 = 0, & \text{in $X$}, \\
   v = 0, & \text{on $M$}, \\
   B_1^3v+\frac{1}{3}H=0, & \text{on $M$} .
  \end{cases}
 \end{equation}
\end{defn}

Equivalently, $\hg$ is the unique metric such that $\hg\rv_{TM}=h$, $\hQ_4=0$, and $\hH=0$.  It is clear that $\mT_3$ depends only on $(M^3,h)$ and the conformal fill-in $(X^4,[g])$, and not a specific choice of metric in $X$.  There are two ways in which $\mT_3^3$ can be regarded as a critical $Q$-curvature.  The first is its conformal transformation formula.

\begin{prop}
 \label{prop:mT_conf}
 Let $(X^4,g)$ be a compact Riemannian manifold with boundary $(M^3,h)=(\partial X,g\rv_{TM})$ and suppose that $\ker L_4\cap\ker B_0^3\cap\ker B_1^3=\{0\}$.  Given $w\in C^\infty(M)$, set $\hh=e^{2w}h$.  Then
 \[ e^{3w}\hmT_3^3 = \mT_3^3 + \mB_3^3w . \]
\end{prop}

\begin{proof}
 Without loss of generality, we may suppose that the metric $g$ is such that $Q_4=0$ and $H=0$.  Let $\hg=e^{2v}g$ be the unique metric such that $\hg\rv_{TM}=\hh$, $\hQ_4=0$, and $\hH=0$.  It follows from the conformal transformation formula for $Q_4$ and Proposition~\ref{prop:l4_conf_crit} that $L_4v=0$, $B_0^3v=w$, and $B_1^3v=0$.  Therefore
 \[ e^{3w}\hmT_3 = \frac{1}{2}e^{3w}\hT_3^3 = \frac{1}{2}\left(\hT_3^3 + B_3^3v\right) = \mT_3^3 + \mB_3^3w . \qedhere \]
\end{proof}

Proposition~\ref{prop:mB_conf} and Proposition~\ref{prop:mT_conf} together imply that the total $\mT$-curvature is a conformal invariant of $(M^3,h)$ as the boundary of $(X^4,[g])$.  Indeed, the total $\mT$-curvature is closely related to the Euler characteristic of $X$ (cf.\ \cite{Anderson2001b,ChangQingYang2006}).

\begin{prop}
 \label{prop:mT_gauss_bonnet}
 Let $(X^4,g)$ be a compact Riemannian manifold with boundary $M^3=\partial X$ and suppose that $\ker L_4\cap\ker B_0^3\cap\ker B_1^3=\{0\}$.  Then
 \[ 8\pi^2\chi(X) = \int_X \lv W\rv^2 + 2\oint_M \mT_3^3 - \frac{2}{3}\oint_M \tr A_0^3 . \]
\end{prop}

\begin{proof}
 The Gauss--Bonnet--Chern formula states that (cf.\ \cite{ChangQing1997a} and Lemma~\ref{lem:cf_changqing_local})
 \[ 8\pi^2\chi(X) = \int_X \left( \lv W\rv^2 + Q_4\right) + \oint_M \left(T_3^3 - \frac{2}{3}\tr A_0^3\right) . \]
 Choosing $g$ such that $Q_4=0$ and $H=0$ yields the result.
\end{proof}

Proposition~\ref{prop:mT_conf} and Proposition~\ref{prop:mT_gauss_bonnet} state that the $\mT$-curvature has the same properties as the fractional $Q$-curvature $Q_3$ on the boundary of a four-dimensional Poincar\'e--Einstein manifold.  In light of Theorem~\ref{thm:asymptotics}, it is unsurprising that these curvatures are the same (cf.\ \cite{ChangQingYang2006}).

\begin{prop}
 \label{prop:mT_is_Q}
 Let $(X^4,M^3,g_+)$ be a Poincar\'e--Einstein manifold with $2\not\in\sigma_{pp}(-\Delta_{g_+})$ and let $h$ be a representative of the conformal boundary.  Then $\mT_3^3=Q_3$.
\end{prop}

\begin{proof}
 Let $r$ be the geodesic defining function associated to $h$.  Let $v$ be the unique solution to $-\Delta_{g_+}v=3$ with
 \[ v = \log r + A + Br^3 \]
 for $A,B\in C^\infty(X)$ even in $r$ and $A\rv_M=0$.  Since $(L_4)_{g_+}=\Delta_{g_+}^2+2\Delta_{g_+}$, it also holds that $(L_4)_{g_+}v+6=0$.  By definition~\cite{FeffermanGraham2002}, $Q_3=3B\rv_M$.

 Consider now the compactified metric $g=r^2g_+$.  It follows from the conformal transformation formulae for $Q_4$ and $L_4$, the evenness of $A$ in $r$, and the fact that $M$ is minimal with respect to $g$ that $v-\log r$ solves~\eqref{eqn:mT_ext}.  Note that $\ker L_4\cap\ker B_0^3\cap\ker B_1^3=\{0\}$ since $2\not\in\sigma_{pp}(-\Delta_{g_+})$.  Hence, by Proposition~\ref{prop:l4_conf_crit},
 \[ \mT_3^3 = \frac{1}{2}\left(T_3^3 + B_3^3(v-\log r)\right) . \]
 By direct computation, $T_3^3=0$ and $B_3^3(v-\log r)=6B\rv_M$.
\end{proof}

In the setting of four-manifolds with boundary, the Paneitz operator and its boundary operators are related to the functionals $\mF,\mG\colon C^\infty\to\bR$ defined by
\begin{align*}
 \mF(u) & = \int_X \left(u\,L_4u + 2Q_4u\right) + \oint_M \left( fB_3^3u + 2T_3^3f + \psi B_2^3u + 2T_2^3\psi\right), \\
 \mG(u) & = \oint_M \left(B_1^3(u)B_2^3(u) + T_2^3B_1^3(u) + T_1^3B_2^3(u)\right) .
\end{align*}
for $f=B_0^3u$ and $\psi=B_1^3u$.  These constitute two of the integrals which naturally arise in log determinant formulas in this setting.  More precisely, Chang and Qing showed~\cite{ChangQing1997a} that the log determinant can be written as a linear combination of eleven functionals, one of which is $b_2+\frac{1}{12}D$, where
\begin{align*}
 b_2(u) & = \frac{1}{4}\int_X \left(u\,L_4u + 2Q_4u\right) + \frac{1}{2}\oint_M \left(fP_3u + 2Tf\right), \\
 D(u) & = \oint_M P_3^2u
\end{align*}
for $u\in C^\infty(X)$, where $f=B_0^3u$.  These functionals are defined precisely as they appear in~\cite[Theorem~3.3]{ChangQing1997a}.  As an immediate consequence of Lemma~\ref{lem:cf_changqing_local}, we find that $b_2+\frac{1}{12}D$ differs from a particular combination of $\mF$ and $\mG$ by integration against pointwise conformal invariants.

\begin{cor}
 \label{cor:cf_changqing}
 Let $(X^4,g)$ be a compact Riemannian manifold with boundary $M^3=\partial X$.  Then
 \[ 4b_2(u) + \frac{1}{3}D(u) = \mF(u) - \mG(u) + \oint_M \left[ \frac{1}{4}\lv A_0\rv^2\psi - 8\left(\lp W(\eta,\cdot,\eta,\cdot),A_0 - \frac{2}{3}\tr A_0^3\right)f \right] \]
 for all $u\in C^\infty(X)$, where $f=B_0^3u$ and $\psi=B_1^3u$.
\end{cor}

One reason to regard the $\mF$- and $\mG$-functionals as natural components of the log determinant functional is that they independently satisfy the following ``cocycle'' condition (cf.\ \cite{BransonGover2008b}).

\begin{prop}
 \label{prop:cohomology}
 Let $(X^4,g)$ be a compact Riemannian manifold with boundary $M^3=\partial X$.  Given any $u,v\in C^\infty(X)$ it holds that
 \begin{align*}
  \hmF(v) & = \mF(u+v) - \mF(u), \\
  \hmG(v) & = \mG(u+v) - \mG(u),
 \end{align*}
 where $\hmF$ and $\hmG$ are defined in terms of the metric $\hg=e^{2u}g$.
\end{prop}

\begin{proof}
 This is an immediate consequence of Theorem~\ref{thm:boundary} and Proposition~\ref{prop:l4_conf_crit}.
\end{proof}

The cocycle condition allows one to readily identify the geometric significance of critical points of the underlying functional.

\begin{prop}
 \label{prop:3d_critical_F}
 Let $(X^4,g)$ be a compact Riemannian manifold with boundary $M^3=\partial X$ and set
 \[ [g]_{1,M} = \left\{ u\in C^\infty(X) \suchthat \oint_M e^{3f} = 1, \text{ where $f=B_0^3u$} \right\} . \]
 Suppose that $u\in[g]_{1,M}$ is a critical point of $\mF\colon[g]_{1,M}\to\bR$.  Then the metric $\hg=e^{2u}g$ is such that
 \begin{equation}
  \label{eqn:3d_critical_F}
  \begin{cases}
   \Vol_{\hh}(M) = 1, \\
   \hQ_4 = 0, \\
   \hT_2^3 = 0, \\
   \hT_3^3 = \int Q_4 + \oint T_3^3 .
  \end{cases}
 \end{equation}
\end{prop}

\begin{remark}
 By imposing the additional constraint $B_1^3u=0$, we obtain that the metric $\hg$ associated to a critical point $u$ satisfies $\hT_1^3=0$ instead of $\hT_2^3=0$.
\end{remark}

\begin{proof}
 Note that $[g]_{1,M}$ consists of all functions $u\in C^\infty(X)$ such that $M$ has unit volume with respect to $\hh=e^{2f}h$.  Let $v\in C^\infty(X)$ be such that $\oint B_0^3(v)e^{3f}=0$.  It follows from Proposition~\ref{prop:cohomology} that
 \[ 0 = \frac{1}{2}\left.\frac{d}{dt}\right|_{t=0}\mF(u+tv) = \int_X \hQ_4v + \oint_M \left(\hT_3^3\hB_0^3(v) + \hT_2^3\hB_1^3(v)\right) . \]
 We thus conclude that $\hQ_4=0=\hT_2^3$ and that $\hT_3^3$ is constant.  By Corollary~\ref{cor:boundary_q4_invariant}, this constant must be $\int Q_4+\oint T_3^3$.
\end{proof}

The situation for the $\mG$-functional is more subtle.  Since the $\mG$-functional depends only on the values of $u$ in a neighborhood of the boundary, it is easy to see that its critical points correspond to conformal metrics with vanishing $T_2^3$- and $T_3^3$-curvature.  A more interesting Euler equation is obtained by restricting to Paneitz-harmonic functions on a Poincar\'e--Einstein manifold.

\begin{prop}
 \label{prop:3d_critical_G}
 Let $(X^4,g)$ be a conformal compactification of a Poincar\'e--Einstein manifold $(X^4,M^3,g_+)$.  Suppose that $u\in\ker L_4$ is a critical point of $\mG\colon\ker L_4\to\bR$.  Then the metric $\hg=e^{2u}g$ is such that
 \begin{equation}
  \label{eqn:3d_critical_G}
  2\hmB_1^3\hT_1^3 + \hT_2^3 = 0 .
 \end{equation}
\end{prop}

\begin{remark}
 The conclusion remains true if $(X^4,g)$ has the property that if $u\in\ker L_4\cap\ker B_1^3$, then $B_2^3u=0$.  By Lemma~\ref{lem:mB_uniqueness}, this property ensures that $\mB_1^3$ is well-defined without prescribing $B_0^3u$.  The conclusion is also true for general compact Riemannian four-manifolds with boundary if $\mG$ is restricted to $\ker L_4\cap\ker B_0^3$.
\end{remark}

\begin{proof}
 Let $v\in\ker L_4$.  It follows from Proposition~\ref{prop:cohomology} that
 \[ 0 = \left.\frac{d}{dt}\right|_{t=0}\mG(u+tv) = \oint_M \left(\hT_2^3\hB_1^3(v) + \hT_1^3\hB_2^3(v)\right) . \]
 Set $\hpsi=\hB_1^3(v)$.  By Theorem~\ref{thm:asymptotics}, it holds that $\hB_2^3(v)=2\hmB_1^3(\hpsi)$.  Therefore
 \[ 0 = \oint_M\left(\hT_2^3\hpsi + 2\hT_1^3\hmB_1^3(\hpsi)\right) = \oint_M\left(\hT_2^3 + 2\hmB_1^3\hT_1^3\right)\hpsi . \]
 Since $\hpsi$ is arbitrary, \eqref{eqn:3d_critical_G} holds.
\end{proof}

We now restrict our attention to the $\mF$-functional; with regards to geometric curvature prescription problems, comparison of Proposition~\ref{prop:l4_geometric_euler} and Proposition~\ref{prop:3d_critical_F} shows that the $\mF$-functional should be regarded as the four-dimensional analogue of the volume-normalized energy functional $\mE_4$.  Pursuing this analogy further, one approach to constructing conformal metrics which satisfy~\eqref{eqn:3d_critical_F} is to extremize
\[ \mu_{4,2}(X,M,g) = \inf\left\{ \mF(u) \suchthat u\in C^\infty(X), \oint_M e^{3u}=1 \right\} . \]
This can be regarded as the four-dimensional analogue of $Y_{4,2}(X,M)$, though $\mu(X,M,g)$ is not a conformal invariant: It follows from Proposition~\ref{prop:cohomology} that
\[ \mu_{4,2}(X,M,e^{2v}g) = \mu_{4,2}(X,M,g) - \mF(v) . \]
One can analogously define $\mu_{4,1}(X,M,g)$ by requiring also that $u\in\ker B_1^3$, with similar properties.

An initial step towards extremizing $\mu_{4,2}(X,M,g)$ is to compute the corresponding value and identify the extremal functions in the model space, namely the round upper hemisphere $(S_+^4,d\theta^2)$ (cf.\ \cite{ChangYang1995}).  Like Theorem~\ref{thm:sobolev} in higher dimensions, this result follows from Proposition~\ref{prop:extension} and the sharp fractional Sobolev inequalities on the round spheres~\cite{Beckner1993}.  The corresponding result for $\mu_{4,1}$ was recently proven by Ache and Chang~\cite{AcheChang2015}.

\begin{thm}
 \label{thm:sobolev_critical}
 Let $(S_+^4,d\theta^2)$ be the four-dimensional upper hemisphere with the round metric.  Given any $u\in C^\infty(S_+^4)$, it holds that
 \[ \int_{S_+^4} \left( \left(\Delta u\right)^2 + 2\lv\nabla u\rv^2 \right) \geq 4\oint_{S^3} \psi\oDelta f + \frac{16\pi^2}{3}\log\left(\frac{1}{2\pi^2}\oint_{S^3}e^{f-\of}\right) + \left(4\pi\oint_{S^3}\lv\psi\rv^3\right)^{\frac{2}{3}} , \]
 where $f=B_0^3u$, $\psi=B_1^3$, and $\of=\frac{1}{2\pi^2}\oint f$.  Moreover, equality holds if and only if $L_4u=0$, $f(x)=a_1+\ln(1+x\cdot\xi_1)$ and $\psi(x)=a_2(1+x\cdot\xi_2)^{-1}$ for constants $a_1,a_2\in\bR$ and points $\xi_1,\xi_2\in B_1(0)$.
\end{thm}

\begin{proof}
 By Proposition~\ref{prop:extension}, it holds that
 \[ \mE_4(u) \geq 2\oint_{S_3} \left(f\,P_3f + \psi\,P_1\psi\right) \]
 with equality if and only if $L_4u=0$.  On the other hand, Beckner~\cite{Beckner1993} showed that
 \begin{align*}
  \oint_{S^3} f\,P_3f & \geq \frac{8\pi^2}{3}\log\left(\frac{1}{2\pi^2}\oint_{S^3} e^{f-\of}\right), \\
  \oint_{S^3} \psi\,P_1\psi & \geq \left(2\pi^2\right)^{1/3}\left(\oint_{S^3} \lv\psi\rv^3\right)^{\frac{2}{3}}
 \end{align*}
 with equality if and only if $f(x)=a_1+\ln(1+x\cdot\xi_1)$ and $\psi(x)=a_2(1+x\cdot\xi_2)^{-1}$ for constants $a_1,a_2\in\bR$ and points $\xi_1,\xi_2\in B_1(0)$. Combining these inequalities with the formula for the energy $\mE_4$ on the round upper hemisphere yields the result.
\end{proof}

In fact, the extremizing functions $u\in C^\infty(S_+^4)$ can be written down explicitly, and they correspond to conformal factors for the adapted metrics~\cite{CaseChang2013} associated to the Einstein metrics on $S^3$; for details, see~\cite{AcheChang2015}.

\bibliographystyle{abbrv}
\bibliography{../bib}
\end{document}